\documentclass[11pt]{article}
\usepackage{a4wide}
\usepackage{url}
\usepackage{enumitem}
\usepackage{booktabs}
\usepackage{caption}
\usepackage{multirow}
\usepackage{amsmath, amssymb}
\usepackage{amsthm}
\usepackage{subfigure}
\usepackage{algorithm}
\usepackage{algorithmicx,algorithm, algpseudocode}
\usepackage{pst-all}
\allowdisplaybreaks
\usepackage[affil-it]{authblk}
\usepackage{lipsum}
\usepackage[numbers]{natbib}

\usepackage{tikz}
\usetikzlibrary{calc, positioning}

\def\rev#1{{#1}}

\newtheorem{theorem}{Theorem}
\newtheorem{conjecture}{Conjecture}

\newtheorem{definition}{Definition}

\def\arcset{{\cal A}}
\def\nodeset{{\cal V}}
\def\pathset{{\cal P}}
\def\network{{\cal N}}
\def\statespace{{\cal S}}
\def\relaxedstatespace{{\cal G}}

\makeatletter
\def\blfootnote{\xdef\@thefnmark{}\@footnotetext}
\makeatother

\begin{document}
\pagestyle{plain}

\title{\Large \bf Arc Flow Formulations Based on Dynamic Programming: Theoretical Foundations and Applications}
{\author{Vin\'icius L. de Lima$^{(1),*}$, Cl\'audio Alves$^{(2)}$, Fran\c{c}ois Clautiaux$^{(3)}$, Manuel Iori$^{(4)}$, Jos\'e M. Val\'erio de Carvalho$^{(2)}$}
\affil{(1) Instituto de Computa{\c{c}}{\~a}o, Universidade Estadual de Campinas (Brazil)\\
(2) Escola de Engenharia / Centro Algoritmi, Universidade do Minho (Portugal) \\
(3) Universit{\'e} de Bordeaux, IMB UMR CNRS 5251, Inria Bordeaux Sud-Ouest (France)\\
(4) DISMI, Universit{\`a} di Modena e Reggio Emilia (Italy) \\
}

\blfootnote{\noindent $*$ corresponding author.}
\blfootnote{e-mail addresses: v.loti@ic.unicamp.br; claudio@dps.uminho.pt; francois.clautiaux@math.u-bordeaux.fr; manuel.iori@unimore.it; vc@dps.uminho.pt.} 

\date{}
\maketitle
\vspace*{-5ex}
\noindent

\begin{abstract}
Network flow formulations are among the most successful tools to solve optimization problems. Such formulations correspond to determining an optimal flow in a network. One particular class of network flow formulations is the arc flow, where variables represent flows on individual arcs of the network. For {\cal NP}-hard problems, polynomial-sized arc flow models typically provide weak linear relaxations and may have too much symmetry to be efficient in practice. Instead, arc flow models with a pseudo-polynomial size usually provide strong relaxations and are efficient in practice. The interest in pseudo-polynomial arc flow formulations has grown considerably in the last twenty years, in which they have been used to solve many open instances of hard problems. A remarkable advantage of pseudo-polynomial arc flow models is the possibility to solve practical-sized instances directly by a Mixed Integer Linear Programming solver, avoiding the implementation of complex methods based on column generation.

In this survey, we present theoretical foundations of pseudo-polynomial arc flow formulations, by showing a relation between their network and Dynamic Programming (DP). This relation allows a better understanding of the strength of these formulations, through a link with models obtained by Dantzig-Wolfe decomposition. The relation with DP also allows a new perspective to relate state-space relaxation methods for DP with arc flow models. We also present a dual point of view to contrast the linear relaxation of arc flow models with that of models based on paths and cycles. To conclude, we review the main solution methods and applications of arc flow models based on DP in several domains such as cutting, packing, scheduling, and routing.
\end{abstract}


\section{Introduction}

Optimization problems can assume different characterizations, each allowing the reduction of the original problem to other optimization problems, leading to, possibly, different solution methods.
One of such characterizations is based on a network, that is, a directed graph with costs on the arcs, and is the base of the well-known network flow problems (see, e.g., \citet{AMO93}).
A {\em network flow problem} is an optimization problem that requires to determine an optimal flow on a network, by satisfying flow conservation on each node and possible additional side constraints on the flow on the arcs.
The list of real-life problems that can be solved as network flow problems is extensive, including not only direct applications (i.e., the network is an input of the problem), as vehicle routing, telecommunication network planning, and train scheduling (see, e.g., \citet{GLR19}, \citet{M06}, and \citet{CT12}), but also indirect ones (i.e., the network must be constructed), as cutting and packing, scheduling, and project management (see, e.g., \citet{DI19}, \citet{KIL19}, and \citet{RJMR20}).

In general, any formulation that corresponds to solving a network flow problem is called {\em network flow formulation}.
Following \citet{AMO93}, the two main classes of network flow formulations are the path (and cycles) flow formulations and the arc flow formulations. {\em Path (and cycles) flow formulations} (called path flow formulations for short in the following) have variables corresponding to flow on paths and cycles of the network, and are frequently associated with set-covering, set-packing, or set-partitioning models. In contrast, {\em arc flow formulations} have variables corresponding to flow on individual arcs of the network.
The Flow Decomposition Theorem by \citet{AMO93} guarantees that the two formulations produce equivalent models when based on the same network, in the sense that any solution of one model can be mapped into a solution of the other model.

Arc flow formulations have a linear number of variables with respect to the number of arcs in the network.
Often, such formulations can be used to solve medium-sized instances directly by a commercial Mixed Integer Linear Programming (MILP) solver, avoiding the implementation of complex methods. 
On the other hand, path flow formulations may have a much larger number of variables, as the number of paths and cycles can be exponential with respect to the number of arcs of the network. Such exponential path flow formulations are typically solved by sophisticated methods based on column generation and branch-and-price algorithms (see, e.g., \citet{DDS06} and \citet{SV13}).

Many $\cal{NP}$-hard problems can be formulated as compact (i.e., polynomial-sized) models. Although their size is favorable, these models are usually associated with weak linear relaxations and very symmetrical solution spaces, leading to an overly extensive enumeration on branch-and-bound algorithms.
To overcome such inefficiency, one may rely on the Dantzig-Wolfe (DW) decomposition (see \citet{DW61}), which can lead to models with stronger linear relaxations and less symmetry.
The models resulting from DW decomposition usually have exponentially more variables than the original model and are solved by column generation.
When the associated pricing problem is solved by {\em Dynamic Programming} (DP), the model can be seen as a path flow formulation, where each variable corresponds to a path in the network inherent from the DP problem. From this observation, one can devise an equivalent arc flow model based on the DP network. Originated by a DW decomposition, the resulting arc flow model will, possibly, have a strong linear relaxation and a less symmetrical solution space, compared to the original polynomial-sized model.

To obtain strong models by a DW decomposition, one may have to pay the price of having an increased complexity on the pricing problem, as such complexity usually becomes pseudo-polynomial or exponential. When the size of the DP network from the pricing problem is pseudo-polynomial, the resulting arc flow model can still be used in practice to solve medium-sized instances by means of a MILP solver. On the other hand, when the DP network is too large (possibly exponential), one can obtain smaller (yet still strong) pseudo-polynomial arc flow models relying on a state-space relaxation of the DP. Different state-space relaxations lead to arc flow models with different sizes and different strength, allowing a flexible possibility to balance size and strength of arc flow models.
 
\subsubsection*{Seminal works and previous surveys}
To the best of our knowledge, \citet{FF58a} were the first to propose a pseudo-polynomial arc flow model, which was based on a time-expanded network and was used to solve the maximal dynamic flow problem. Already in the sixties, \citet{S68} studied the relation between DP and Integer Programming, by modeling the knapsack problem as a pseudo-polynomial network flow problem. More than ten years later, \citet{W77} proposed a general methodology to build packing and covering models from the set of feasible solutions of DP problems. As an application, he introduced a pseudo-polynomial arc flow model based on a DP network of the knapsack problem for the cutting stock problem (CSP), a strongly ${\cal NP}$-hard problem. 
Due to the context, \citet{W77} did not present computational experiments for the model, and such formulations did not get much attention for decades. It was more than twenty years later that \citet{V99} independently proposed this model, used it to solve to (integer) optimality open instances in the CSP literature, and showed that it is equivalent to the Gilmore-Gomory model~\cite{GG61,GG63}, a path flow model obtained from a DW decomposition. 

Since the work by \citet{V99}, the popularity of pseudo-polynomial arc flow formulations based on DP has grown considerably, and several different applications to ${\cal NP}$-hard problems have appeared in the literature. Large instances have been solved for the first time to proven optimality by arc flow models in several cutting, packing, and scheduling problems. In addition, arc flow formulations have effectively modeled complex problems involving, for instance, multiple-stage cutting processes or the presence of setups.

Some previous surveys considered arc flow models on specific applications (see, e.g., \citet{V02} and \citet{DIM16} for bin packing and cutting stock problems, and \citet{CAF16} for vehicle routing problems with multiple trips).

\citet{S09} presented a survey on problems based on dynamic flow, which is also called "flows over time". One of the main formulations for this class of problems is the arc flow model based on time-expanded networks, introduced by \citet{FF58a}, which is derived from DP graphs of pseudo-polynomial size. A generalization of such networks is the layered graph approach which considers, for instance, capacity-indexed networks for vehicle routing problems. Modeling techniques and efficient solution methods for pseudo-polynomial arc flow formulations based on layered graphs were surveyed by \citet{GLR19}. 
Other surveys, such as that by \citet{SV13}, focused instead on how arc flow formulations can be used to devise effective branch-and-price algorithms. \citet{CCZ10} and \citet{LS14} studied techniques to transform large extended formulations into compact ones. In the case of network flow formulations, this occurs when path flow models are transformed into equivalent arc flow models with (possibly exponentially) fewer variables.

\subsubsection*{Contents}
In this survey, we extend the previous studies to provide a base for some reasons behind the (primal) strength of pseudo-polynomial arc flow formulations. Based on the Flow Decomposition Theorem, we show how DW decomposition can derive arc flow models with a strong linear relaxation. We also discuss how the state-space relaxation method may allow one to obtain a balance between strength and complexity when constructing arc flow models. This relation between state-space relaxation and arc flow formulations can implicitly relate different arc flow models in the literature (see Section \ref{sec:state-space-relaxation}). We provide a dual insight that may explain arc flow models' practical computational efficiency over their equivalent path flow models, which is the richer description of the dual space. The relevance is that a better description of the dual space leads to less primal degeneracy, which is headwind for Linear Programming (LP) simplex (vertex) algorithms used in branch-and-bound techniques. We also discuss the main solution methods to solve large-scale arc flow models and present the main applications studied in the literature.

The remainder of this paper is organized as follows. Section \ref{sec:foundations} presents the theoretical foundations of network flow formulations and DP. The relationship between arc flow models with a strong relaxation and DW decomposition is presented in Section \ref{sec:dantzig_wolfe_decomposition}. In Section \ref{sec:state-space-relaxation}, we present a discussion on how the state-space relaxation can be used to obtain smaller arc flow models. In Section \ref{sec:dual_insight}, we show how arc flow formulations present a more detailed dual space than path flow formulations, which leads to a better convergence of the LP solution. The state-of-the-art methods to solve large-scale arc flow models and the main applications are shown in Sections \ref{sec:methods} and \ref{sec:applications}, respectively. Finally, our concluding thoughts are presented in Section \ref{sec:conclusion}.

\section{Network Flow Formulations and Dynamic Programming} \label{sec:foundations}

A {\em network} $\network$ is composed of a directed graph $G = (\nodeset,\arcset)$, where $\nodeset$ is the set of nodes (vertices) and $\arcset \subseteq \nodeset \times \nodeset$ is the set of arcs, and of a cost function $c \colon \arcset \to \mathbb{R}$ that maps each arc $(u,v)$ to a cost $c(u,v)$. In addition, a network has two special nodes in $\nodeset$, called {\em source} ($v_{source}$) and {\em sink} ($v_{sink}$), such that no arcs in $\arcset$ enters $v_{source}$ or leaves $v_{sink}$. Depending on the context, we may denote the set of nodes and the set of arcs of a network $\network$ as $\nodeset(\network)$ and $\arcset(\network)$, respectively.

A {\em flow} in a network is a function $f \colon \arcset \to \mathbb{R}_+$ that satisfies the flow conservation constraints:
\begin{align}
&\sum_{ (u,v) \in \arcset} f((u,v)) = \sum_{(v,w) \in \arcset} f((v,w)), & \forall v \in \nodeset \setminus \{v_{source}, v_{sink} \}, 
\end{align}
i.e., the flow entering a node is equal to the flow leaving it, except for $v_{source}$ and $v_{sink}$, where there is only outgoing and incoming flow, respectively.
 By the flow conservation constraints, the flow leaving $v_{source}$ is equal to the flow entering $v_{sink}$.

\begin{definition}
A network flow problem is an optimization problem which requires to determine a flow $f$ that satisfies side constraints and optimizes $\sum_{(u,v) \in \arcset} c(u,v)f((u,v))$.
In addition, any formulation that corresponds to solving a network flow problem is referred to as a network flow formulation.
\end{definition}

A {\em path} in a graph is a sequence of arcs which joins a sequence of nodes. A {\em cycle} is a path in which the first and the last nodes are the same. We denote by $\pathset$ the set of all paths from $v_{source}$ to $v_{sink}$, and by $\arcset(p)$ the set of arcs of a path $p \in \pathset$. The cost of a path $p \in \pathset$ is given by $\tilde{c}_p = \sum_{(u,v) \in \arcset(p)} c{(u,v)}$. In practice, the solution of network flow problems can be given as a set of paths and cycles with a positive flow or as the flow on each arc of the network. Following \citet{AMO93}, this gives rise to two different classes of network flow formulations: path flow formulations and arc flow formulations.
\begin{definition}
Path flow formulations are network flow formulations where decision variables correspond to the non-negative flow on each path and each cycle of the network.
\end{definition}

Path flow formulations are also referred to in the literature as {\em path-based formulations}.
In path flow formulations, the flow conservation is implicitly imposed on the variables. Since $|\pathset(\network)|$ can be exponential with respect to $|\arcset(\network)|$, path flow formulations usually have a huge number of variables. Nonetheless, such formulations have been successfully solved in the literature by column generation based algorithms (see, e.g., \citet{PU14}). 

\begin{definition}
Arc flow formulations are network flow formulations where decision variables correspond to the non-negative flow on each arc of the network.
\end{definition}

In contrast to path flow formulations, arc flow formulations impose the flow conservation explicitly as linear constraints. For instance, let variable $\varphi_{(u,v)} \in \mathbb{R}_+$, for each arc $(u,v) \in \arcset$, represent the flow on arc $(u,v)$, and let variable $z \in \mathbb{R}_+$ represent the total flow on the network. Flow conservation constraints can be formulated as:
\begin{align}
\label{eq:network_flow_conservation}
\sum_{ (u,v) \in \arcset} \varphi_{(u,v)} - \sum_{ (v,w) \in \arcset } \varphi_{(v,w)} & = \begin{cases} -z, &\text{ if } v = v_{source}, \\ z, &\text{ if } v = v_{sink} , \\ 0, &\text{ otherwise}, \end{cases} & \forall v \in \nodeset.
\end{align}

The underlying matrix from \eqref{eq:network_flow_conservation} is totally unimodular (see, e.g., \citet{NW88}), which is an interesting property, as LP models with a totally unimodular constraint matrix have the {\em integrality property}, i.e., every vertex of the corresponding polytope is integer, implying the existence of an integer optimal solution if the right-hand side of the constraints is integer.

Although flow conservation constraints \eqref{eq:network_flow_conservation} must be explicitly included, differently from path flow formulations, arc flow formulations have a polynomial number of variables with respect to network size. An important result that relates these two classes of formulations is the following:

\begin{theorem}[Flow Decomposition Theorem (\citet{AMO93})] \label{thm:flow_decomposition}
Every path and cycle flow has a unique representation as non-negative arc flows. Conversely, every non-negative arc flow can be represented as a path and cycle flow (though not necessarily uniquely) with the following two properties:
\begin{enumerate}[label=(\alph*)]
\item Every directed path with a positive flow connects the source to the sink;
\item At most $|\nodeset|$ + $|\arcset|$ paths and cycles have nonzero flow; out of these, at most $|\arcset|$ cycles have nonzero flow.
\end{enumerate}
\end{theorem}

Theorem \ref{thm:flow_decomposition} proves the equivalence between arc flow and path flow formulations: an arc flow solution can be transformed into positive flow on a set of paths and cycles, and a solution of path flow formulations can be decomposed into flow on individual arcs. Consequently, arc flow and path flow formulations based on the same network produce the same linear relaxation bounds for integer problems.

A fundamental network flow problem is the {\em longest path problem} (LPP): find a longest path, i.e., a unitary flow of maximum cost, from $v_{source}$ to $v_{sink}$. When the network does not have cycles of negative cost, a longest path can be found in polynomial time with respect to the network size. In particular, when the network is acyclic, a longest path can be found in $O(|\arcset|)$ by a topological ordering of the nodes (see, e.g., \citet{AMO93}).
The LPP often appears in solution methods for more complex network flow problems. Such problems have, for instance, side constraints characterized by generalized upper and lower bounds on the flow on subsets of arcs (see, e.g., \citet{CHMA17}).
In the remainder of this section, we formally define Dynamic Programming and its relation to acyclic network flow formulations and the LPP.

\subsection{Dynamic Programming and Arc Flow Formulations}
(Discrete) Dynamic Programming (DP) is a well-known method, proposed in the fifties (see, e.g., \citet{B57}), to solve combinatorial optimization problems that can be decomposed into a sequence of decisions, often represented by {\em stages}, each corresponding to a decision step. DP models are defined by a state space $\statespace$, where each state $s \in \statespace$ is characterized by a set of entities. Each state is defined, for instance, by the subset of clients already visited in routing problems (see, e.g., \citet{DDS92}) or by the level of usage of a given resource in, e.g., cutting, packing or scheduling problems (see, e.g., \citet{CH95}).

A key aspect in DP modeling is the fact that the description of a state $s \in \statespace$ should be sufficient to recognize the admissible decisions for $s$, which is coined as the {\em no-memory property of DP}. Then, when stages are considered, problems are divided into a sequence of sub-problems, so that the solution of a sub-problem depends only on the sub-problems from the previous stage. Nowadays, some authors changed the perspective of DP from the stage-dependent description to a table-filling method, as described by, e.g., \citet{AMO93}.

A DP model can be represented by a recursive function, in which the value of a state is only based on the values of the precedent states and the cost (contribution) of the decisions leading to that state. Formally, for each state $s \in \statespace$, let $\Delta(s)$ be the set of precedent states, and $c{(r,s)} \in \mathbb{R}$ be the cost required to move from $r\in \Delta(s)$ to $s$. The state space $\statespace$ contains two special states, $s_0$ and $ s^{*}$, representing the initial condition of the recursion and the optimal solution of the problem, respectively. No state precedes $s_0$ ($\Delta(s_0) = \emptyset$) and $s^*$ does not precede any other state ($s^{*} \notin \Delta(s), \forall s \in \statespace$). The DP recursion we consider is given by:
\begin{align}
\label{eq:dp_recursion}
& f_{DP}(s) = \begin{cases} \max_{\{r \in \Delta(s)\}} \{ f_{DP}(r) + c{(r,s)} \}, & \text{ if } s\neq s_0, \\ 0, &\text{ if } s = s_0, \end{cases} & \forall s \in \statespace.
\end{align}

A priori, the DP recursion \eqref{eq:dp_recursion} denotes a maximization problem, but it can easily be considered as a minimization problem by inverting the sign of all decision costs.
Based on recursion \eqref{eq:dp_recursion}, a DP model can be defined on a directed acyclic graph (DAG), where vertices and arcs correspond to states and decisions, respectively. In this characterization, the DP solution is given by a longest path, that is, a path from $s_0$ to $s^*$ with maximum total cost (profit). In some cases, the description of the recursion produces states that are not in any solution from the initial state to the goal state, and such states can be disregarded.

The previous characterization of DP is a classical definition, and one of the most used by the integer programming literature. Nevertheless, other definitions of DP were proposed. For instance, \citet{MRC90} defined DP over a directed acyclic hypergraph, to obtain a better polyhedral characterization of a class of combinatorial optimization problems. This generalization was later used to solve network flow problems (see, e.g., \citet{CSVV18}).

Following our DP definition, the DAG from a DP model produces the {\em dynamic programming network} $\network_{DP}$. Each node is associated with a state, that is, $\nodeset(\network_{DP}) \equiv \statespace$, where states $s_0$ and $s^*$ correspond to nodes $v_{source}$ and $v_{sink}$, respectively. The set of arcs is defined by $\arcset(\network_{DP}) = \{ (r, s) \mid s \in \statespace, r \in \Delta(s) \}$, and the cost $c{(r,s)}$ of an arc $(r,s) \in \arcset$ is equivalent to the decision cost $c{(r,s)}$ for moving from state $r$ to state $s$.

The fact that DP can be seen as a LPP in the DP network provides a generic recipe to transform a DP model into the arc flow model given by:
\begin{align}
\label{eq:DP_arcflow_of}
\max & \sum_{(u,v) \in \arcset(\network_{DP})} c{(u,v)}\varphi_{(u,v)}, \\
\label{eq:DP_arcflow_flow_conservation}
\text{s.t.:} & \sum_{(u,v) \in \arcset(\network_{DP})} \varphi_{(u,v)} - \sum_{(v,w) \in \arcset(\network_{DP})} \varphi_{(v,w)} = \begin{cases} -1, &\text{ if } v = v_{source}, \\ 1, &\text{ if } v = v_{sink} , \\ 0, &\text{ otherwise}, \end{cases} & \forall v \in \nodeset(\network_{DP}), \\
\label{eq:DP_arcflow_flow_domain}
& \varphi_{(u,v)} \in \{0, 1\}, & \forall (u,v) \in \arcset(\network_{DP}).
\end{align}

Model \eqref{eq:DP_arcflow_of}--\eqref{eq:DP_arcflow_flow_domain} considers flow conservation constraints with a unitary flow ($z = 1$) and side constraints impose that the flow on each arc is binary. The objective function is to maximize the cost of the selected path, given by the sum of the contributions of the state decisions. As this model has only flow conservation constraints, it has the integrality property, because, as previously discussed, its underlying matrix is totally unimodular. This property is linked with the fact that the LPP on a DAG can be solved in linear time with respect to the number of arcs (see, e.g., \citet{B57}).

In strongly ${\cal NP}$-hard problems, DP often provides pseudo-polynomial algorithms to solve sub-problems that determine partial plans. Examples of such sub-problems are: finding a single cutting/packing pattern on cutting and packing problems (see, e.g., \citet{V02}); obtaining a schedule in a single machine in scheduling problems (see, e.g., \citet{KIL19}); determining a route of a single vehicle on vehicle routing problems (see, e.g., \citet{PU14}).

The DP network from the sub-problems (partial plans) may be used to derive arc flow models to solve the main problem (global plan). If the global plan is a strongly ${\cal NP}$-hard problem (which is often the case), these models do not have the integrality property, unless ${\cal P} = {\cal NP}$, but they usually provide strong relaxations. We provide examples of DP networks of two problems that are often associated with partial plans.

\subsection{Example on the Knapsack Problem} \label{sec:dp_example_knapsack}

Let us consider an example on the {\em knapsack problem} (KP): given a knapsack with capacity $W \in \mathbb{Z}_+$ and a set $I = \{1,2,\ldots,n\}$ of items, each item $i \in I$ with a weight $w_i \in \mathbb{Z}_+$, a profit $p_i \in \mathbb{Z}_+$, and a maximum number of copies $d_i \in \mathbb{Z}_+$, determine the number $x_i$ (with $0 \leq x_i \leq d_i$) of copies of each item $i \in I$ in the solution, so that $\sum_{i \in I} w_i x_i \leq W$ and $\sum_{i \in I} p_i x_i$ is maximum (see, e.g., \citet{MT90}). The optimal solution value $f_{KP}(i,W')$ of a knapsack sub-problem that considers the items $\{1,2,\ldots,i\}$ and a partial capacity $W'$ can be recursively computed by:
\begin{align}
\label{eq:knapsack_recursion_1}
& f_{KP}(i, W') = \begin{cases} \max_{ l \in \{0,1,\dots,\min\{d_i,\lfloor W'/w_i \rfloor\} \}} \{ f_{KP}(i-1,W'- w_i l) + p_i l \}, & \text{ if } i > 0 \text{ and } W' > 0, \\ 0, & \text{ otherwise}. \end{cases}
\end{align} 

The optimal KP solution value is given by $OPT(I, W) = \max_{\{ 0 < W'\leq W\}}f_{KP}(|I|, W')$ and it is associated with a DP model where the state space $\statespace$ is a subset of $\{ (i,W') \mid i \in I, W' \in \{0,1,\ldots,W \} \} \cup \{ s^* \}$, each partial set of items corresponds to a stage, the initial condition state is $s_0 = (0,0)$, and the goal is represented by a dummy state $s^*$, so that the precedent states in $\Delta(s^*)$ are related to the recursion of $OPT(I,W)$. The set of precedent states of each state $(i,W')\in\statespace \setminus \{ s^*\}$, where $i>0$, is $\Delta((i,W')) = \{ (i-1,W'-w_il) \mid w_il \leq W', l \in \{0,1,\ldots,d_i\} \}$, and the decision cost for reaching $(i, W')$ from a precedent state $(i-1, W'-w_il) \in \Delta((i,W'))$ is $c{((i-1, W'-w_il),(i,W'))} = p_il$. The set of precedent states of $s^*$ is $\Delta(s^*) = \{ (|I|,W') \mid W' \in \{0,1,\ldots,W \} \}$ and the cost to move from a precedent state $(|I|, W') \in \Delta(s^*)$ to $s^*$ is $c{((|I|,W'), s^*)} = 0$.

Let ${\cal N}_{KP}$ be a DP network for the KP. The set of nodes $\nodeset({\network}_{KP}) \equiv \statespace$ corresponds to the set of states, where $v_{source}$ and $v_{sink}$ are associated with $s_0$ and $s^*$, respectively. The set of arcs is $\arcset({\cal N}_{KP}) = \cup_{\{i\in I, l \in \{ 0,1,\dots,\min\{d_i,\lfloor W'/w_i \rfloor \}\} } \arcset_{il} \cup \arcset_S$. For each $i\in I$, the set $\arcset_{il} = \{ ((i-1, W'-w_il), (i, W')) \mid (i, W') \in \statespace, (i-1, W'-w_il) \in \Delta((i, W'))\} $ contains the arcs corresponding to the decision of choosing $l$ copies of $i$ in the solution, and these arcs have profit $p_il$. The set $\arcset_S = \{ ((n, W'), v_{sink}) \mid (n,W') \in \statespace \}$ contains arcs, called {\em loss arcs}, that link the nodes from the last stage to $v_{sink}$, and their profit is $0$.

As an example, Figure \ref{fig:knapsack_dp} shows the network for a KP with capacity $W=8$, and three items having $w_1 = 4$, $w_2 = 3$, $w_3 = 2$, $p_1 = 8$, $p_2 = 5$, $p_3 = 4$, and $d_1 = d_2 = d_3 = 1$. There is a node $s^*$ representing $v_{sink}$, a node $(i,W')$ for each state, and $v_{source} = (0,0)$. Arcs $\arcset_S$ are depicted as dotted links, and for every item $i\in I$, the arcs in $\arcset_{i1}$ and $\arcset_{i0}$ are depicted as full and dashed links, respectively. The arcs that do not belong to any path from $(0,0)$ to $s^*$ are disregarded. By definition, the only non-zero profit arcs are the ones from $\arcset_{i1}$, and their profits are shown in the figure. The bold arcs are from the longest path, i.e., the optimal solution, which contains items $1$ and $2$ and has total profit $13$.

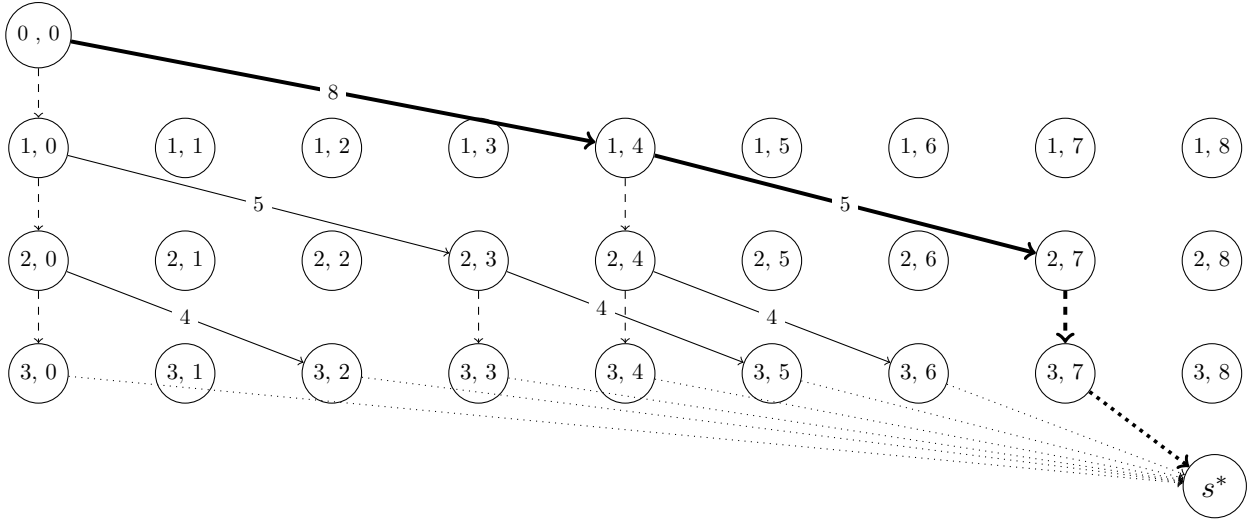
\begin{figure}[!hbt]
\center
\begin{tikzpicture}[darkstyle/.style={circle,draw,fill=gray!40,minimum size=20}, scale = 0.75, transform shape]

  \foreach \x in {0,...,8}
    \foreach \y in {1,...,3} 
       {\pgfmathtruncatemacro{\ypos}{5 - \y}
       \node [scale = 1.5, darkstyle, fill=white] (\y\x) at (2.6*\x,2*\ypos) {};}

 \foreach \x in {0,...,8}
   \foreach \y in {1,...,3} 
      {
      \node [scale = 1.05] at (\y\x) {\y , \x};} 

  \node [scale = 1, darkstyle, fill=white] (00) at (0,10) {0 , 0}; 

  \node[scale = 1.4, darkstyle, fill=white] (goal) at (20.85,2) {$s^*$};
	
	\draw[dotted, ->] (30) edge (goal);
	\draw[dotted, ->] (32) edge (goal);
	\draw[dotted, ->] (33) edge (goal);
	\draw[dotted, ->] (34) edge (goal);
	\draw[dotted, ->] (35) edge (goal);
	\draw[dotted, ->] (36) edge (goal);
	\draw[dotted, ->, line width = 1.4] (37) edge (goal);

	\draw[->,line width=1.5] (00) -- (14) node[midway, fill=white]{8};
	\draw[dashed, ->] (00) -- (10);

	\draw[dashed, ->] (10) -- (20);
	\draw[->] (10) -- (23) node[midway,fill=white]{5};

	\draw[dashed, ->] (14) -- (24);
	\draw[->, line width = 1.5] (14) -- (27) node[midway,fill=white]{5};

	\draw[dashed, ->] (20) -- (30);
	\draw[->] (20) -- (32) node[midway,fill=white]{4};

	\draw[dashed, ->] (24) -- (34);
	\draw[->] (24) -- (36) node[midway,fill=white]{4};

	\draw[dashed, ->] (23) -- (33);
	\draw[->] (23) -- (35) node[pos = 0.4,fill=white]{4};

	\draw[dashed, ->, line width = 1.5] (27) -- (37);

\end{tikzpicture}
\caption{Example of network $\network_{KP}$ for the knapsack problem.}
\label{fig:knapsack_dp}
\end{figure}

The network $\network_{DP}$ can be used in the generic recipe previously discussed, to derive an arc flow model for the KP with $O(|I|W)$ nodes and arcs. Although the resulting arc flow model has pseudo-polynomial size and has the integrality property, solving it with LP algorithms is usually not as fast as solving the problem by DP. However, the LP model can be a base for interesting results. For instance, \citet{AB92} extended this model to solve the knapsack separation problem, which is a generic tool to generate cutting planes for MILP formulations.

\subsection{Example on the Elementary Shortest Path Problem with Resource Constraints} \label{sec:example_ESPPRC}

\begin{figure}[!hbt]
\center
\begin{tikzpicture}[darkstyle/.style={circle,draw,fill=gray!40,minimum size=20}, scale = 0.75, transform shape]

  \node [scale = 1.0, darkstyle, fill=white] (n1) at (0,0) {\{\}, 0};
  \node [scale = 1.0, darkstyle, fill=white] (n2) at (-9,-2.0) {\{1\}, 1};
  \node [scale = 1.0, darkstyle, fill=white] (n3) at (-3,-2.0) {\{2\}, 2};
  \node [scale = 1.0, darkstyle, fill=white] (n4) at (3,-2.0) {\{3\}, 3};
  \node [scale = 1.0, darkstyle, fill=white] (n5) at (9,-2.0) {\{4\}, 4};
  \node [scale = 1.0, darkstyle, fill=white] (n6) at (-11,-5) {\{1,2\}, 2};
  \node [scale = 1.0, darkstyle, fill=white] (n7) at (-9,-5) {\{1,3\}, 3};
  \node [scale = 1.0, darkstyle, fill=white] (n8) at (-7,-5) {\{1,4\}, 4};
  \node [scale = 1.0, darkstyle, fill=white] (n9) at (-4.5,-5) {\{1,2\}, 1};
  \node [scale = 1.0, darkstyle, fill=white] (n10) at (-1.5,-5) {\{2,3\}, 3};
  \node [scale = 1.0, darkstyle, fill=white] (n12) at (1.5,-5) {\{1,3\}, 1};
  \node [scale = 1.0, darkstyle, fill=white] (n13) at (4.5,-5) {\{2,3\}, 2};
  \node [scale = 1.0, darkstyle, fill=white] (n14) at (9,-5) {\{1,4\}, 1};
  \node [scale = 1.0, darkstyle, fill=white] (n15) at (-6,-8.5) {\{1,2,3\}, 3};
  \node [scale = 1.0, darkstyle, fill=white] (n16) at (0,-8.5) {\{1,2,3\}, 2};
  \node [scale = 1.0, darkstyle, fill=white] (n17) at (6,-8.5) {\{1,2,3\}, 1};
  \node [scale = 1.4, darkstyle, fill=white] (goal) at (0,-11) {$s^*$};

  \draw[->] (n1) -- (n2) node[midway, fill=white]{-5};
  \draw[->] (n1) -- (n3) node[midway, fill=white]{-3};
  \draw[->, line width = 1.5] (n1) -- (n4) node[midway, fill=white]{1};
  \draw[->] (n1) -- (n5) node[midway, fill=white]{-2};
  \draw[->] (n2) -- (n6) node[midway, fill=white]{6};
  \draw[->] (n2) -- (n7) node[midway, fill=white]{2};
  \draw[->] (n2) -- (n8) node[midway, fill=white]{-2};
  \draw[->] (n3) -- (n9) node[midway, fill=white]{-7};
  \draw[->] (n3) -- (n10) node[midway, fill=white]{-1};
  \draw[->] (n4) -- (n12) node[midway, fill=white]{2};
  \draw[->, line width = 1.5] (n4) -- (n13) node[midway, fill=white]{-5};
  \draw[->] (n5) -- (n14) node[midway, fill=white]{-2};
  \draw[->] (n6) -- (n15) node[midway, fill=white]{-1};
  \draw[->] (n7) -- (n16) node[midway, fill=white]{-5};
  \draw[->] (n9) -- (n15) node[pos=0.6, fill=white]{2};
  \draw[->] (n10) -- (n17) node[midway, fill=white]{2};
  \draw[->] (n12) -- (n16) node[midway, fill=white]{6};
  \draw[->, line width = 1.5] (n13) -- (n17) node[midway, fill=white]{-7};

  \draw[dotted,->] (n1) -- (goal);
  \draw[dotted,->] (n2) -- (goal);
  \draw[dotted,->] (n3) -- (goal);
  \draw[dotted,->] (n4) -- (goal);
  \draw[dotted,->] (n5) -- (goal);
  \draw[dotted,->] (n6) -- (goal);
  \draw[dotted,->] (n7) -- (goal);
  \draw[dotted,->] (n8) -- (goal);
  \draw[dotted,->] (n9) -- (goal);
  \draw[dotted,->] (n10) -- (goal);
  \draw[dotted,->] (n12) -- (goal);
  \draw[dotted,->] (n13) -- (goal);
  \draw[dotted,->] (n14) -- (goal);
  \draw[dotted,->] (n15) -- (goal);
  \draw[dotted,->] (n16) -- (goal);
  \draw[dotted,->, line width = 1.5] (n17) -- (goal);

\end{tikzpicture}
\caption{\rev{Example of network $\network_{ESPPRC}$ for the elementary shortest path problem with resource constraints.}}
\label{fig:ESPPRC_DP}
\end{figure}
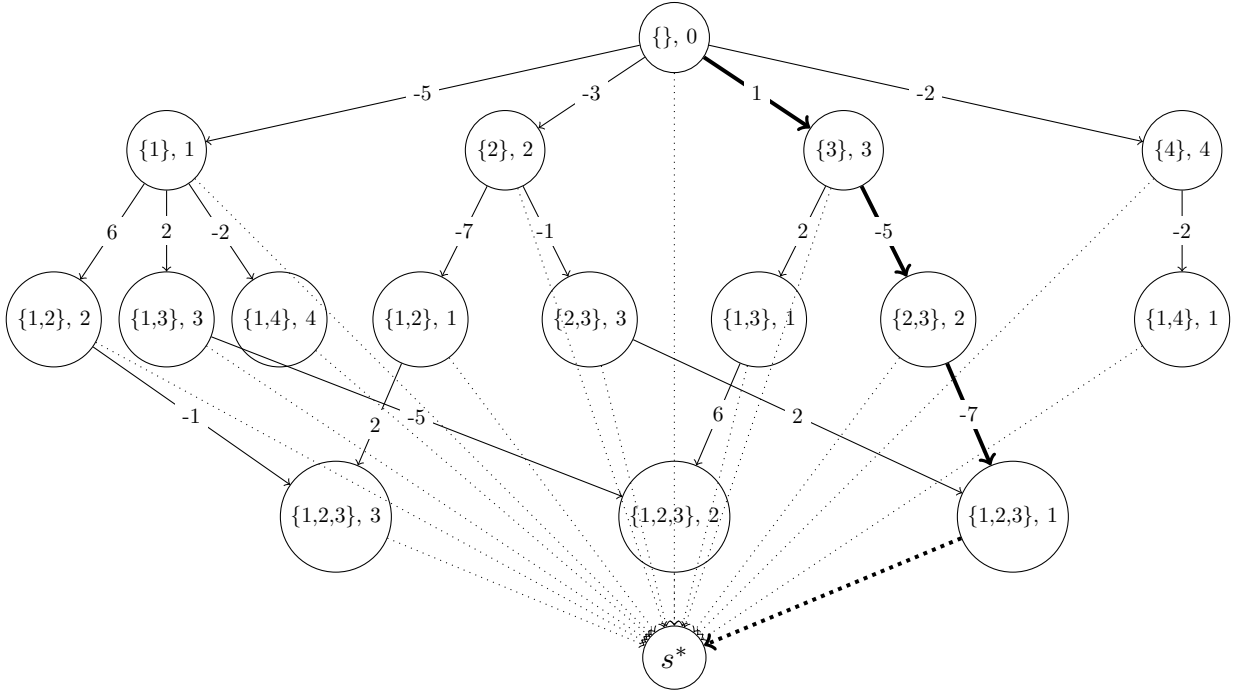

The {\em elementary shortest path problem with resource constraints} (ESPPRC) is linked with determining the best route of a single vehicle in vehicle routing problems (see, e.g., \citet{ID05}). In the ESPPRC there is a resource capacity and a directed graph. Each vertex is associated with a resource consumption and each arc is associated with a cost. The objective is to find a path with minimum cost such that: (i) each vertex is visited at most once (elementary constraint); and (ii) the total resource consumption from the vertices in the path does not exceed the resource capacity (resource constraint).

In the context of vehicle routing, the vertices are given by $I = \{0,1,\ldots,n\}$, where $0$ represents the depot and the other $n$ vertices represent the clients. The solution must start and finish at the depot. Generally, the ESPPRC might consider a set of resources, but in this example we are concerned only with a single resource $W$, which is associated with the load capacity of a vehicle. Let $e_{(i,j)} \in \mathbb{R}$ be the cost of arc $(i,j) \in I \times I$ and $w_i \in \mathbb{R}_+$ be the resource consumption of $i\in I$.

Given a set of clients $S \subseteq I \setminus \{ 0 \}$ that satisfies the resource constraint (i.e., $\sum_{j \in S} w_j \leq W$), we denote by $f_{ESPPRC}(S, i)$ the cost of a path of minimum cost from the depot $0$ to a client $i \in S$ that visits every client in $S$ exactly once. A client $j \in S \setminus \{ i \}$ that precedes $i$ in an optimal path related to $f_{ESPPRC}(S, i)$ is one that minimizes the sum of the cost $e_{(j,i)}$ from $j$ to $i$ and the cost $f_{ESPPRC}(S \setminus \{i\}, j)$ of a path of minimum cost that visits every client in $S \setminus \{ i\}$ exactly once and finishes in $j \in S \setminus \{ i\}$. This relation can be mapped into the following recursion:
\begin{align}
\label{eq:recursion_cvrp_espprc}
f_{ESPPRC}(S, i)& = \begin{cases} \min_{j \in S \setminus \{ i \}} \{ f_{ESPPRC}(S \setminus \{ i\}, j) + e_{(j,i)} \}, & \text{ if } S \neq \emptyset \text{ and } i > 0, \\ 0, & \text{ otherwise}. \end{cases}
\end{align}

The optimal ESPPRC solution value is given by $ OPT(I, e, r) = \min \{ f_{ESPPRC}(S, i) + e_{(i,0)} \mid S \subseteq I, i\in S, \sum_{j \in S} w_j \leq W \} $. The recursion can be transformed into a maximization problem, by inverting the sign of each $e_{(i,j)}$, resulting in a DP model. The state space $\statespace$ is a subset of $\{ (S, i) \mid S \subseteq I, i \in S, \sum_{j \in S} w_j \leq W \} \cup \{ s^* \}$, where the initial state is $s_0 = (\emptyset, 0)$ and the goal state $s^*$ has precedent states related to the recursion of $OPT(I,e,r)$. For each state $(S, i) \in \statespace \setminus \{ s^*\}$, the set of precedent states is $\Delta((S, i)) = \{ (S \setminus \{i\}, j) \mid j \in S \setminus \{i\}\}$, and the decision cost to move from state $(S \setminus \{ i \}, j) \in \Delta((S, i))$ to $(S, i)$ is $c{((S \setminus \{i\}, j), (S, i))} = e_{(j,i)}$. The set of precedent states of $s^*$ is $\Delta(s^*) = \statespace \setminus \{ s^* \}$, and the cost to move from state $(S, i) \in \Delta(s^*)$ to $s^*$ is $c{((S,i),s^*)} = e_{(i,0)}$.

Let $\network_{ESPPRC}$ be the DP network for the ESPPRC obtained from the state space $\statespace$ presented above. The set of nodes $\nodeset(\network_{ESPPRC}) \equiv \statespace$ corresponds to the state space, where $v_{source}$ and $v_{sink}$ are associated with $(\emptyset, 0)$ and $OPT(I,e,r)$, respectively. The set of arcs $\arcset(\network_{ESPPRC}) = \cup_{i\in I} \arcset_i $ contains arcs $\arcset_i = \{ ((S \setminus \{ i \}, j), (S,i)) \in \nodeset \times \nodeset \}$ related to the decision of visiting client $i \in I \setminus \{ 0 \}$, and arcs $\arcset_0 = \{((S, j), v_{sink}) \in \nodeset \times \nodeset\} $ related to the decision of returning to the depot.

\rev{Figure \ref{fig:ESPPRC_DP} illustrates the network for an example with $W = 6$ and $4$ clients, with $w_1 = 1$, $w_2 = 2$, $w_3 = 2$, and $w_4 = 5$. In this example, the costs are the following: $c_{(0,1)}=-5$, $c_{(0,2)}=-3$, $c_{(0,3)}=1$, $c_{(0,4)}=-2$, $c_{(1,0)}=0$, $c_{(1,2)}=6$, $c_{(1,3)}=2$, $c_{(1,4)}=-3$, $c_{(2,0)}=0$, $c_{(2,1)}=-7$, $c_{(2,3)}=-1$, $c_{(2,4)}=5$, $c_{(3,0)}=0$, $c_{(3,1)}=2$, $c_{(3,2)}=-5$, $c_{(3,4)}=5$, $c_{(4,0)}=0$, $c_{(4,1)}=-2$, $c_{(4,2)}=5$, $c_{(4,3)}=5$. On each arc, we present its associated cost, except for the arcs that enter $s^*$, which for the sake of conciseness we consider their cost to be $0$ and present such arcs as dotted lines. The arcs of the optimal solution are highlighted in bold, and they correspond to visit the sequence of clients $3, 2$, and $1$, and returning to the depot. The optimal solution has value $-11$.
}

The network $\network_{ESPPRC}$ can be the base of an arc flow model for the ESPPRC, following the generic recipe previously presented. However, the number of subsets $S \subseteq I$ such that $\sum_{j\in S}w_j \leq W$ is exponential, implying that the state space presented for the ESPPRC is also exponential. Consequently, the corresponding arc flow model has an exponential number of variables and constraints, and it is not practical to solve it directly. Practical techniques that can be used to overcome this inconvenience are described in Section \ref{sec:state-space-relaxation}. \rev{In particular, Section \ref{sec:state_space_relaxation_CVRP} presents a network of pseudo-polynomial size for a version of the ESPPRC where the elementary constraints are partially relaxed.}

\section{Dantzig-Wolfe Decomposition and Network Flow Formulations} \label{sec:dantzig_wolfe_decomposition}

Several challenging problems admit compact MILP models that, although having a reasonable (polynomial) size, are usually associated with weak linear relaxations. In such cases, reformulation methods can obtain models with stronger linear relaxation. One of such methods is the well-known {\em Dantzig-Wolfe} (DW) {\em decomposition} (see \citet{DW61}), which has been successfully used in many applications. The resulting models may be stronger, but they usually have an exponential number of variables, and complex methods are required to solve them in practice. In the following, we show how models obtained from DW decomposition can be the base to derive equivalent arc flow models with smaller (and possibly practical) size.

First, we show how the DW decomposition can be applied to Linear Programming (LP) models. Let us express the feasible solution region of the LP model as:
\begin{align}
\label{eq:dw_lp_region}
X_{LP} = \{ x \in \mathbb{R}^n_+ \mid Ax \geq b, x \in X \},
\end{align}

\noindent where $A \in \mathbb{R}^{m \times n}$, $b \in \mathbb{R}^{m}$, and $X$ is a bounded convex polytope of a set of constraints. According to the Minkowski's Theorem, any point in X can be represented as a convex combination of the vertices (extreme points) of X (for the general case of an unbounded polytope, it is also necessary to consider a non-negative linear combination of the extreme rays of the polytope, see, e.g., \citet{NW88}). Then, by defining ${\cal Q}(X)$ as the set of vertices of $X$, it follows:
\begin{align}
\label{eq:minkowvvski_reformulation}
X = \{ x \in \mathbb{R}_+^n \mid x = \sum_{q \in {\cal Q}(X)} q\lambda_q, \sum_{q \in {\cal Q}(X)} \lambda_q = 1, \lambda_q \geq 0, \forall q \in {\cal Q}(X) \}.
\end{align}

By substituting $x$ as in \eqref{eq:minkowvvski_reformulation} in \eqref{eq:dw_lp_region}, the feasible region in the space of the variables of the reformulated model becomes:
\begin{align}
X_{LPDW} = \{ \lambda \in \mathbb{R}_+^{|{\cal Q}(X)|} \mid \sum_{q \in {\cal Q}(X)} (Aq)\lambda_q \geq b, \sum_{q \in {\cal Q}(X)} \lambda_q = 1 \}.
\end{align}

The reformulated problem, called {\em DW model}, has a variable associated with each vertex of $X$, and a new set of constraints, usually referred to as convexity constraints (see, e.g., \citet{LD05}). \rev{The DW model is just another way of expressing the same solution space. Therefore,} the decomposition applied to LP models preserves the optimal solution value. 

However, in the context of MILP, to obtain linear models stronger than the original linear relaxation, the integrality constraints must be implicitly considered in the reformulated variables. A possibility is to \rev{impose the integrality constraints just in the set X, reformulating} over the vertices of $Conv\{ x\in X \text{ and integer} \}$, the convex hull of the integer points in $X$ (see, e.g., \citet{NW88}), to optimize over the set:

\begin{equation}
X_{DW} = \{ x \in \mathbb{R}_{+}^{n} | A x \geq b, x \in Conv \{ x \in X \mbox{ and integer} \}\}. \label{XDW}
\end{equation}

Note indeed that $X_{DW} \subseteq X_{LP}$, and the relation can be strict when the set~$X$ does not have the integrality property. \rev{When that happens, the reformulated model is stronger. The strength of a model} can have a strong impact on the search for the optimal integer solution. Typically the length of the search is smaller with stronger models, leading to smaller computational times.
For further details, the reader is referred to \citet{NW88}, which discusses at length the importance of deriving stronger models in the context of MILP.

\rev{
A typical issue from the DW decomposition is that the resulting model may have an exponential number of variables. This issue is usually addressed by relying on the column generation method, a technique to solve LP models with a large number of variables (see, e.g., \citet{LD05}). The method solves a restricted version of the LP model with only a subset of variables by iteratively solving a {\em pricing problem}. The pricing generates non-basic variables (columns) that are candidates to improve the current restricted model. If no such variable exists, then the current basis of the restricted model is optimal for the original model, and the method halts.
Column generation algorithms are complex, but it pays off to solve models resulting from a DW decomposition, because they may be stronger \rev{when the set $X$ does not have the integrality property and integrality is enforced in the pricing problem.}
}

When the pricing problem of a DW model can be solved by DP, each column of the model can be associated with a path in the DP network, and the model can be seen as a path flow formulation. In such cases, the coefficients of a column are given as contributions from the arcs of the corresponding path. In addition, due to the no-memory property of the DP, for each column, the contribution $\ell_{(u,v)}^k \in \mathbb{R}$ of an arc $(u,v)$ to a row coefficient $k \in {\cal C}$, where ${\cal C}$ is the set of constraints, is independent of other arcs. Then, given a DP network ${\cal N}$, and $d_k \in \mathbb{R}$, for every constraint $k \in {\cal C}$, we obtain a general path flow formulation:
\begin{align}
\label{eq:general_path_formulation_of}
\min & \sum_{p\in {\cal P(N)}} \tilde{c}_p\lambda_p, & \\
\label{eq:general_path_formulation_constraints}
\text{s.t.:} & \sum_{p \in {\pathset(\network)}} \sum_{(u,v) \in \arcset(p)} \ell_{(u,v)}^k \lambda_p \geq d_k, & \forall k \in {\cal C}, \\
\label{eq:general_path_formulation_domain}
& \lambda_p \in \mathbb{Z}_+, & \forall p \in {\cal P(\network)}.
\end{align}

The objective function \eqref{eq:general_path_formulation_of} minimizes the total cost from the paths in the solution, and \eqref{eq:general_path_formulation_constraints} is a set of general linear constraints. As discussed previously, due to the flow decomposition theorem, the path flow formulation \eqref{eq:general_path_formulation_of}--\eqref{eq:general_path_formulation_domain} can be reformulated as the following arc flow formulation:
\begin{align}
\label{eq:general_arc_flow_formulation_of}
\min ~ & \sum_{(u,v) \in {\arcset(\network)}} c{(u,v)}\varphi_{(u,v)}, & \\
\label{eq:general_arc_flow_conservation_constraints}
\text{s.t.: } & \sum_{ (u,v) \in {\arcset(\network)}} \varphi_{(u,v)} - \sum_{ (v,w) \in {\arcset(\network)} } \varphi_{(v,w)} = \begin{cases} -z, &\text{ if } v = v_{source}, \\ z, &\text{ if } v = v_{sink} , \\ 0, &\text{ otherwise}, \end{cases} & \forall v \in \nodeset(\network), \\
\label{eq:general_arc_flow_general_constraints}
& \sum_{(u,v) \in {\arcset(\network)}} \ell_{(u,v)}^k \varphi_{(u,v)} \geq d_k, & \forall k \in {\cal C}, \\
\label{eq:general_arc_flow_formulation_domain1}
& z \geq 0, & \\ 
\label{eq:general_arc_flow_formulation_domain2}
& \varphi_{(u,v)} \in \mathbb{Z}_+, & \forall (u,v) \in {\arcset(\network)}.
\end{align}

The objective function \eqref{eq:general_arc_flow_formulation_of} minimizes the total cost from the arcs in the solution, \eqref{eq:general_arc_flow_conservation_constraints} are the flow conservation constraints, and \eqref{eq:general_arc_flow_general_constraints} are general linear constraints, adapted from \eqref{eq:general_path_formulation_constraints}.

DW decomposition can generate path flow formulations that are usually associated with strong linear relaxations, and equivalent arc flow formulations can be derived from the DP network of the underlying pricing problem. As pointed before, to obtain stronger relaxations, one should reformulate over the convex hull of the integer points in a polytope that does not have the integrality property. Although a stronger linear relaxation is obtained, this kind of reformulation leads to pricing problems that are ${\cal NP}$-hard.

On the other hand, for arc flow models, the price to pay to obtain a stronger relaxation (i.e., to obtain a set equivalent to $X_{DW}$) is to use a set of constraints that defines a DP network with a pseudo-polynomial or an exponential number of arcs. 

Usually, arc flow models of practical size can be derived when the resulting pricing problem has pseudo-polynomial complexity. However, when the pricing problem is strongly ${\cal NP}$-hard, the resulting arc flow model is generally too large to be solved in practice.
When the pricing problem from path flow formulations is too difficult in practice, one may rely on relaxation methods. We show in Section \ref{sec:state-space-relaxation} one of such relaxation methods, which can be used to derive smaller networks and obtain arc flow formulations that are still strong and have practical size. In the next sections, we present two examples of arc flow formulations derived from DW decomposition, with pseudo-polynomial and exponential size, respectively.

\subsection{Example on the Cutting Stock Problem} \label{sec:example1_csp}

We present an example of an arc flow model derived from a DW decomposition for the {\em Cutting Stock Problem} (CSP). In the CSP, a roll of width $W$ and a set $I=\{1,2,\ldots,n\}$ of items are given, such that each item $i \in I$ has a positive demand $d_i$ and a width $w_i$. The objective is to cut the demand of all items from the minimum number of rolls.

Given an upper bound $K$ on the minimum number of rolls, following \citet{MT90}, the CSP can be formulated as: 
\begin{align}
\label{eq:kantorovich_of}
\min & \sum _{j=1}^K y_j, & \\
\label{eq:kantorovich_knapsack}
\text{s.t.: } & \sum_{i=1}^n w_i x_{ij} \leq W y_j, & j = 1,2,\ldots,K, \\
\label{eq:kantorovich_demand}
& \sum_{j=1}^K x_{ij} \geq d_i, & \forall i \in I, \\
\label{eq:kantorovich_domain1}
& y_j \in \{ 0, 1 \}, & j = 1,2,\ldots,K, \\
\label{eq:kantorovich_domain2}
& x_{ij} \in \mathbb{Z}_+, & \forall i \in I, j=1,2,...,K.
\end{align}

Each roll $j=1,2,\ldots,K$ is associated with a binary variable $y_j$ that is equal to $1$ if and only if roll $j$ is used in the solution. The integer variable $x_{ij}$ represent the number of copies of item $i\in I$ that is cut from roll $j$. The objective function \eqref{eq:kantorovich_of} minimizes the number of rolls cut in the solution. Constraints \eqref{eq:kantorovich_knapsack} ensure that the size of each roll is satisfied and that no item is cut from an unused roll, whereas constraints \eqref{eq:kantorovich_demand} ensure that the demand of each item is satisfied.

According to \citet{MT90}, the optimal solution value of the linear relaxation of \eqref{eq:kantorovich_of}--\eqref{eq:kantorovich_domain2} is equivalent to the continuous lower bound that is obtained from the minimum length required to cut all of the demand from a stock with unlimited width, i.e., $\sum_{i \in I} w_id_i / W$. This bound is known to be weak in practice, and its worst-case performance ratio asymptotically tends to $1/2$. To obtain a stronger bound, a DW decomposition can be applied to the model above (see, e.g., \citet{V98}).

Constraints \eqref{eq:kantorovich_knapsack} correspond to $K$ identical fractional knapsack polytopes. When integrality constraints \eqref{eq:kantorovich_domain1} and \eqref{eq:kantorovich_domain2} are taken into account, the $K$ identical polytopes correspond to $K$ identical sets $P_{CSP}$ of integer knapsack solutions. A DW decomposition that includes constraints \eqref{eq:kantorovich_knapsack}, \eqref{eq:kantorovich_domain1}, and \eqref{eq:kantorovich_domain2} in the sub-problem results in the following set-covering model:
\begin{align}
\label{eq:csp_set_covering_of}
\min & \sum_{p \in P_{CSP}} \lambda_{p}, & \\
\label{eq:csp_set_covering_demand}
\text{s.t.: } & \sum_{p \in P_{CSP}} a_{ip} \lambda_{p} \geq d_i, & \forall i \in I, \\
\label{eq:csp_set_covering_domain}
& \lambda_{p} \in \mathbb{Z}_+, & \forall p \in P_{CSP}.
\end{align}

In the CSP, the integer knapsack solutions from $P_{CSP}$ are called {\em cutting patterns} (or just patterns, for short, in the following), and each of them represents the cutting of a single piece of stock. The objective function \eqref{eq:csp_set_covering_of} minimizes the number of cutting patterns. Constraints \eqref{eq:csp_set_covering_demand} ensure that the demand of each item is satisfied, where $a_{ip}$ is the number of copies of item $i$ in pattern $p$.

In cutting and packing problems, a {\em proper pattern} is a pattern that respects the maximum number of copies of the items, whereas a {\em non-proper pattern} does not. The set-covering model by \citet{GG61, GG63} for the CSP is similar to \eqref{eq:csp_set_covering_of}--\eqref{eq:csp_set_covering_domain}, but it includes non-proper patterns that are obtained from the polytope of an {\em unbounded knapsack problem} (UKP), a variant of the KP where each item has an unlimited number of copies. The linear relaxation of \eqref{eq:csp_set_covering_of}--\eqref{eq:csp_set_covering_domain}, often referred to as {\em proper relaxation} (see, e.g., \citet{KRSK15}), is stronger than the linear relaxation of the model by \citet{GG61, GG63}, but the optimal solution values of the corresponding MILP models are equal.

The linear relaxation of \eqref{eq:csp_set_covering_of}--\eqref{eq:csp_set_covering_domain} is strong and often the number of rolls in the optimal solution is the rounded up optimal solution value from this relaxation. In fact, there is a conjecture related to the strength of this relaxation (see, e.g., \citet{CDDIR15} and \citet{KRSK15}):

\begin{conjecture}[Modified Integer Round-Up Property (MIRUP)]
The difference between the optimal solution value of the CSP and the rounded-up solution value of the linear relaxation of \eqref{eq:csp_set_covering_of}--\eqref{eq:csp_set_covering_domain} is at most one.
\end{conjecture}

The pricing problem of the set-covering model is a KP (where each item has $d_i$ copies), which, as shown in Section \ref{sec:dp_example_knapsack}, can be solved by DP, implying on the existence of a DP network ${\cal N}_{KP}$. As every column from the set-covering model can be represented as a path in ${\cal N}_{KP}$, this network can be the base of an arc flow model for the CSP:
\begin{align}
\label{eq:csp_cambazard_of}
\min & ~ z, \\
\label{eq:csp_cambazard_flow_conservation}
\text{s.t.:} & \sum_{(u,v) \in \arcset({\cal N}_{KP})} \varphi_{(u,v)} - \sum_{(v,w) \in \arcset({\cal N}_{KP})} \varphi_{(v,w)} = \begin{cases} -z, &\text{ if } v = v_{source}, \\ z, &\text{ if } v = v_{sink} , \\ 0, &\text{ otherwise}, \end{cases} & \forall v \in \nodeset(\network_{KP}), \\
\label{eq:csp_cambazard_demand}
& \sum_{(u,v) \in \arcset_{i}({\cal N}_{KP})} \ell_{(u,v)}^i \varphi_{(u,v)} \geq d_i, & \forall i \in I, \\
& z \in \mathbb{Z}_+, & \\
\label{eq:csp_cambazard_domain}
& \varphi_{(u,v)} \in \mathbb{Z}_+, & \forall (u,v) \in \arcset({\cal N}_{KP}).
\end{align}

The objective function \eqref{eq:csp_cambazard_of} minimizes the total flow. Constraints \eqref{eq:csp_cambazard_flow_conservation} impose the flow conservation and constraints \eqref{eq:csp_cambazard_demand} are related to the demand of each item. Following the definitions in Section \ref{sec:dp_example_knapsack}, we further define $\arcset_i(\network_{KP}) = \cup_{l \in \{0,1,\ldots, d_i\}} \arcset_{il}(\network_{KP})$, for each $i \in I$, and the contribution $\ell_{(u,v)}^i$ of arc $(u,v) \in \arcset_i(\network_{DP})$ is the number of copies of $i$ associated with $(u,v)$. The arc flow model \eqref{eq:csp_cambazard_of}--\eqref{eq:csp_cambazard_domain}, corresponds to the arc flow model proposed by \citet{CO10}, under the name {\em DP-flow}, for the {\em bin packing problem} (BPP), a particular case of the CSP where $d_i = 1$ for every $i \in I$. The DP-flow is equivalent to the set-covering model \eqref{eq:csp_set_covering_of}--\eqref{eq:csp_set_covering_domain} (see, e.g., \citet{DI19}), and it follows the proper relaxation.

In this example, we started from a model that is associated with a weak linear relaxation. A DW decomposition resulted in a model with a stronger relaxation, but with the drawback of having a potentially exponential number of variables. Then, we showed how the DP network related to the pricing problem of the exponential model derives an equivalent pseudo-polynomial network flow model. Practical successful applications of this idea are shown in Section \ref{sec:applications}. 

\subsection{Example on the Capacitated Vehicle Routing Problem} \label{sec:example_cvrp}

In the {\em capacitated vehicle routing problem} (CVRP), we are given a set $K$ of identical vehicles and a set $I=\{0,1,\ldots,n\}$ of vertices where vertex $0$ corresponds to a depot, and vertices $1,\ldots, n$ correspond to $n$ clients. Each vehicle has a load capacity $W$, each client $i\in I \setminus \{ 0 \}$ has a non-negative demand $w_i$, and each pair of vertices $i,j \in I$ is associated with a cost $c{(i,j)}$. The CVRP aims at determining a routing plan with the minimum total cost, such that: (i) exactly $|K|$ routes are considered; (ii) each route starts and ends at the depot; (iii) each client is visited exactly once; and (iv) the total demand from the clients of a route does not exceed the load capacity $W$.

In this example, we apply a DW decomposition on a compact arc flow model for the CVRP, obtaining a set partitioning model that derives an equivalent arc flow model of exponential size.

Let $I^* = \{ I \setminus \{ 0 \} \} \cup \{ 0^+, 0^- \}$ be the set of clients with two additional vertices $0^+$ and $0^-$, which are copies of the depot, each corresponding to the source and the sink of a route, respectively. Two additional sets of arcs, $\{ (0^+, j) \mid j \in \{I \setminus \{0\}\} \}$ and $\{ (i, 0^-) \mid i \in \{I \setminus \{0\}\} \}$, are considered. The following compact arc flow model, also known as three-index (vehicle-flow) formulation (see, e.g., \citet{ITV14}), solves the CVRP:
\begin{align}
\label{eq:cvrp1_of}
\min & \sum_{k \in K} \sum_{i \in I^*} \sum_{j \in I^* \setminus \{ i \}} c{(i,j)} \phi^k_{(i,j)}, & \\
\label{eq:cvrp1_flow_conservation}
\text{s.t.: } & \sum_{ j \in I^* \setminus \{ i \}} \phi^k_{(i,j)} - \sum_{ j \in I^* \setminus \{ i \} } \phi^k_{(j,i)} = \begin{cases} -1, &\text{ if } i = 0^+, \\ 1, &\text{ if } i = 0^- , \\ 0, &\text{ otherwise}, \end{cases} & \forall k\in K, i \in I^*, \\
\label{eq:cvrp1_client_demand}
& \sum_{k \in K} \sum_{j \in I^* \setminus \{ i \}} \phi^k_{(i,j)} = 1, & \forall i \in I \setminus \{ 0\}, \\
\label{eq:cvrp1_capacity_constraints1}
& \omega_{ik}-\omega_{jk}+W\phi_{(i,j)}^k \leq W - w_j, & \forall i,j \in I^*, i\neq j, k \in K, \\
\label{eq:cvrp1_capacity_constraints2}
& \omega_{ik} \leq W, & \forall i \in I^*, k \in K, \\
& \phi^k_{(i,j)} \in \{ 0,1 \}, & \forall i,j \in I^*, k \in K, \\
\label{eq:cvrp1_domain}
& \omega_{ik} \geq 0, & \forall i \in I^*, k \in K. 
\end{align}

Model \eqref{eq:cvrp1_of}--\eqref{eq:cvrp1_domain} considers a copy of the original graph for each vehicle $k\in K$. The binary variable $\phi_{(i,j)}^k$ is equal to $1$ if and only if vehicle $k$ moves from client $i$ to client $j$, and the variable $\omega_{ik}$ indicates the accumulated demand already distributed by the vehicle $k$ when arriving at client $i$. The objective function \eqref{eq:cvrp1_of} minimizes the total cost of the routes. Constraints \eqref{eq:cvrp1_flow_conservation} guarantee the conservation of a unitary flow (representing a single route) for each of the $|K|$ vehicles. Constraints \eqref{eq:cvrp1_client_demand} guarantee that each client is visited exactly once. Constraints \eqref{eq:cvrp1_capacity_constraints1} and \eqref{eq:cvrp1_capacity_constraints2} model the propagation of the accumulated demand of the load capacity of each vehicle, and they are also used to ensure that the route of each vehicle is a single connected component (path).

Although model \eqref{eq:cvrp1_of}--\eqref{eq:cvrp1_domain} has a polynomial number of variables, it has a large number of symmetries that makes it inefficient in practice (see, e.g., \citet{ITV14}). Each possible route can be attributed to any vehicle, so each solution (routing plan) has $|K|!$ equivalent permutations. Next, we present a DW decomposition that eliminates this symmetry.

Flow conservation constraints \eqref{eq:cvrp1_flow_conservation} correspond to $|K|$ identical polytopes, each containing all paths from $0^+$ to $0^-$, whose vertices are always integer. By including constraints \eqref{eq:cvrp1_capacity_constraints1} and \eqref{eq:cvrp1_capacity_constraints2}, each of the $|K|$ resulting polytopes contains all paths from $0^+$ to $0^-$ satisfying the load capacity, but the integrality property is lost. Let us consider the DW decomposition over the set $P_{CVRP}$ of vertices of the convex hull of constraints \eqref{eq:cvrp1_flow_conservation}, \eqref{eq:cvrp1_capacity_constraints1}, \eqref{eq:cvrp1_capacity_constraints2}, and \eqref{eq:cvrp1_domain}, that represent the set of integer paths that satisfy the load capacity. The resulting DW model is:
\begin{align}
\label{eq:cvrp_set_partitioning_of}
\min & \sum_{p \in P_{CVRP}} \tilde{c}_{p}\lambda_{p}, & \\
\label{eq:cvrp_set_partitioning_demand}
\text{s.t.: } & \sum_{p \in P_{CVRP}} a_{pi} \lambda_{p} = 1, & \forall i \in I \setminus \{ 0 \}, \\
& \sum_{p \in P_{CVRP}} \lambda_{p} = |K|, &\\
\label{eq:cvrp_set_partitioning_domain}
& \lambda_{p} \in \{0,1\}, & \forall p \in P_{CVRP}.
\end{align}

The set partitioning model \eqref{eq:cvrp_set_partitioning_of}--\eqref{eq:cvrp_set_partitioning_domain} provides a strong lower bound and it is among the most successful formulations to solve the CVRP in practice (see, e.g., \citet{PU14} and \citet{PPPU17}). This model associates with each path $p \in P_{CVRP}$ a variable $\lambda_p$. The binary coefficient $a_{ip}$ is equal to $1$ if and only if client $i \in I \setminus \{ 0 \}$ is visited in $p$. The objective function \eqref{eq:cvrp_set_partitioning_of} minimizes the total cost, where the cost $\tilde{c}_p$ of a path $p$ is the sum of the costs of the arcs in $\arcset(p)$. Constraints \eqref{eq:cvrp_set_partitioning_demand} guarantee that each client is visited by exactly one route. As each variable represents a unique path, the symmetry from the model \eqref{eq:cvrp1_of}--\eqref{eq:cvrp1_domain} is eliminated.

The pricing problem associated with model \eqref{eq:cvrp_set_partitioning_of}--\eqref{eq:cvrp_set_partitioning_domain} is equivalent to the ESPPRC. In Section \ref{sec:example_ESPPRC}, we presented the exponential DP network $\network_{ESPPRC}$ for the ESPPRC. Based on this network, and recalling that $\arcset_i(\network_{ESPPRC})$ is the set of arcs that visit $i \in I$, we can reformulate the set partitioning model into the following arc flow model:
\begin{align}
\label{eq:cvrp_exponential_arc_flow_of}
\min & ~ \sum_{(u,v) \in \arcset(\network_{ESPPRC})} c{(u,v)}\varphi_{(u,v)}, \\
\nonumber
\text{s.t.:} & \sum_{(u,v) \in \arcset({\cal N}_{ESPPRC})} \varphi_{(u,v)} - \sum_{(v,w) \in \arcset({\cal N}_{ESPPRC})} \varphi_{(v,w)} & = \begin{cases} -|K|, &\text{ if } v = v_{source}, \\ |K|, &\text{ if } v = v_{sink} , \\ 0, &\text{ otherwise}, \end{cases} \\
\label{eq:cvrp_exponential_arc_flow_conservation}
& & \forall v \in \nodeset(\network_{ESPPRC}), \\
\label{eq:cvrp_exponential_arc_flow_demand}
& \sum_{(u,v) \in \arcset_i({\cal N}_{ESPPRC})} \varphi_{(u,v)} = 1, & \forall i \in I \setminus \{ 0 \}, \\
\label{eq:cvrp_exponential_arc_flow_domain}
& \varphi_{(u,v)} \in \{0,1\}, & \forall (u,v) \in \arcset({\cal N}_{ESPPRC}).
\end{align}

The objective function \eqref{eq:cvrp_exponential_arc_flow_of} minimizes the total cost of the routing plan. Constraints \eqref{eq:cvrp_exponential_arc_flow_conservation} impose the flow conservation, and constraints \eqref{eq:cvrp_exponential_arc_flow_demand} imply that each client is visited exactly once.

As network $\network_{ESPPRC}$ is possibly exponential, the arc flow model \eqref{eq:cvrp_exponential_arc_flow_of}--\eqref{eq:cvrp_exponential_arc_flow_domain} may be too large to be used in practice. Nonetheless, this model is presented so as to provide an example in Section \ref{sec:state_space_relaxation_CVRP} on how a relaxation for the pricing problem of \eqref{eq:cvrp_set_partitioning_of}--\eqref{eq:cvrp_set_partitioning_domain} derives a pseudo-polynomial arc flow model for the CVRP.

\section{State-Space Relaxation on Arc Flow Formulations} \label{sec:state-space-relaxation}

Dynamic programming is often an efficient tool to solve hard problems. However, when the number of states from a DP model is very large, it is not practical to solve it by enumerating all states, and more sophisticated methods are needed. From that observation, \citet{CMT81} proposed a general relaxation procedure for DP, called {\em state-space relaxation}, \rev{which aggregates subsets of states in order to obtain a smaller state space which provides a bound for the original problem. Such relaxation can be embedded in exact solution methods to find the optimal solution for the original state space by (partially) disaggregating the relaxed state space during the search for feasible solutions. Since a state space may allow many different relaxations, it is desirable to determine relaxations that provide a good balance between number of states and strength of the relaxation.}

Formally, a state-space relaxation consists of a mapping function $g \colon \statespace \rightarrow \relaxedstatespace$ between two state spaces, where $|\relaxedstatespace| < |\statespace|$. For every $s \in \statespace$ and $r \in \Delta(s)$, function $g$ must guarantee that $g(r) \in \Delta(g(s))$, and the decision cost of going from $g(r)$ to $g(s)$ is defined as $\bar{c}_{(g(r), g(s))} = \max_{\{p,q \in \statespace\}} \{ c{(p,q)} \mid g(p) = g(r), g(q) = g(s) \}$. The DP recursion of the resulting state-space relaxation is given by: 
\begin{align}
\label{eq:state_space_relaxation_recursion}
& f_{SSR}(g(s)) = \begin{cases} \max_{ \{ p \in \Delta(g(s)) \}} ( f_{SSR}(p) + \bar{c}_{(p, g(s))} ), & \text{ if } s \neq s_0, \\ 0, & \text{ if } s = s_0, \end{cases} & \forall s \in \statespace.
\end{align}

As a relaxation, recursion \eqref{eq:state_space_relaxation_recursion} guarantees that $f_{SSR}(g(s)) \geq f_{DP}(s)$, i.e., it produces an upper bound for the original recursion $f_{DP}$. 
In the context of arc flow formulations induced by DP, the size of the LP formulation is strictly related to the size of the state space. In this case, smaller arc flow formulations can be obtained from state-space relaxations. Solutions obtained with the state-space relaxation may be unfeasible for the original problem. However, there are many cases in which arc flow formulations based on networks from state-space relaxations are guaranteed to produce optimal integer solutions that are feasible for the original problem. The main drawback of state-space relaxations is that they lead to arc flow formulations with weaker linear relaxation. This weakness occurs as the relaxation can profit from paths that are not feasible in the original network, generating, for instance, non-proper patterns in cutting and packing problems or non-elementary routes in vehicle routing problems. \rev{However, in many cases, the reduction on the size of the model pays off the loss in the linear relaxation strength}.

Strong pseudo-polynomial arc flow formulations can be obtained from state-space relaxations of both pseudo-polynomial and exponential state spaces. Motivated by vehicle routing problems, \citet{GLR19} studied modeling and solution methods of a class of pseudo-polynomial arc flow formulations obtained from state-space relaxation of exponential state-spaces, named by the authors as {\em layered graph formulations}.

Next, we present two examples of pseudo-polynomial arc flow formulations obtained from state-space relaxations over a pseudo-polynomial and an exponential state space.

\subsection{Example on the State Space for the Cutting Stock Problem} \label{sec:ssr_csp}

In Section \ref{sec:dp_example_knapsack}, we presented recursion \eqref{eq:knapsack_recursion_1} to solve the KP. This recursion produces network $\network_{KP}$, which has $O(|I|W)$ nodes and $O(|I|W)$ arcs. Then, in Section \ref{sec:example1_csp}, network $\network_{KP}$ was the base of arc flow model \eqref{eq:csp_cambazard_of}--\eqref{eq:csp_cambazard_domain} for the CSP, i.e., the DP-Flow model of \citet{CO10}.

Let $\network_{KP-SSR}$ be the network of a state-space relaxation of $\network_{KP}$ based on the mapping function $g((i,W')) = (W')$, for every state $(i, W') \in \nodeset(\network_{KP})$. This mapping function disregards the dimension related to the partial set of items and merges states representing the same partial capacities, into a single state. Network $\network_{KP-SSR}$ is smaller than $\network_{KP}$, with $O(W)$ nodes, but it may contain paths that violate the maximum number of copies of each item. Then, an arc flow model for the KP, based on $\network_{KP-SSR}$, needs additional side constraints to impose limits on the number of copies of items, while in $\network_{KP}$, such constraints are implicitly imposed by the configuration of the network.

In the CSP, as each demand constraint \eqref{eq:csp_cambazard_demand} imposes only a minimum and not a maximum number of items, there is no problem in considering a network that has paths representing cutting patterns with more copies of a single item than the required one. Hence, the arc flow model obtained by substituting $\network_{KP}$ by $\network_{KP-SSR}$ on model \eqref{eq:csp_cambazard_of}--\eqref{eq:csp_cambazard_domain} still solves the CSP, because unnecessary items copies, if any, can be removed from the solution without affecting its cost. The resulting model is equivalent to the arc flow model for the CSP proposed by \citet{V99}.

The model with $\network_{KP-SSR}$ is smaller than the original model, with $O(|I|+W)$ constraints instead of $O(|I|W)$, but its linear relaxation is weaker. Nonetheless, $\network_{KP-SSR}$ was obtained by \citet{V99} after reducing the DP network of the UKP, implying the graph of $\network_{KP-SSR}$ is a subgraph of the underlying graph from a DP for the UKP. This implies that the linear relaxation of the model with $\network_{KP-SSR}$ is at least as strong as the relaxation of the model by \citet{GG61,GG63}, whose sub-problem is the UKP. Thus, the resulting linear relaxation still follows the MIRUP conjecture, implying that it is still strong. In practice, the arc flow model by \citet{V99}, which is associated with $\network_{KP-SSR}$, is preferable than the model by \citet{CO10}, which is associated with $\network_{KP}$, as its linear relaxation is still strong and it is substantially smaller.

To exemplify, Figure \ref{fig:knapsack_ssr_network} shows the network $\network_{KP-SSR}$ obtained from the example of Figure \ref{fig:knapsack_dp}.
We conclude that state-space relaxation can be considered to reduce the size of pseudo-polynomial arc flow formulations and obtain models with a linear relaxation that is still strong. This modeling technique to derive smaller, but yet efficient, pseudo-polynomial models from the relaxation of pseudo-polynomial state spaces has been mainly used, even without mentioning it, to model one- and two-dimensional cutting and packing problems and scheduling problems (see Section \ref{sec:applications}).

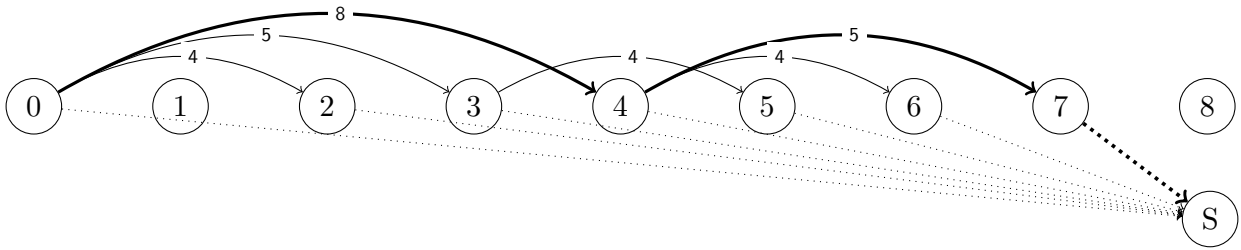
\begin{figure}
\center
\begin{tikzpicture}[darkstyle/.style={circle,draw,fill=gray!40,minimum size=20}, scale = 0.75, transform shape]

  \foreach \x in {0,...,8}{
    \node [scale = 1.4, darkstyle, fill = white] (\x) at (2.6*\x,4) {\x};} 

 \node[scale = 1.4, darkstyle, fill=white] (goal) at (20.85,2) {S};

	\draw[dotted, ->] (0) edge (goal);
	\draw[dotted, ->] (2) edge (goal);
	\draw[dotted, ->] (3) edge (goal);
	\draw[dotted, ->] (4) edge (goal);
	\draw[dotted, ->] (5) edge (goal);
	\draw[dotted, ->] (6) edge (goal);
	\draw[dotted, ->, line width = 1.4] (7) edge (goal);

	\path[every node/.style={font=\sffamily\small}]
    (0) edge[->, bend left, line width = 1.2] node [right, fill = white] {8} (4)
    (0) edge[->, bend left] node [right, fill = white] {5} (3)
    (4) edge[->, bend left, line width = 1.2] node [right, fill = white] {5} (7)
    (0) edge[->, bend left] node [right, fill = white] {4} (2)
    (3) edge[->, bend left] node [right, fill = white] {4} (5)
    (4) edge[->, bend left] node [right, fill = white] {4} (6);

\end{tikzpicture}
\caption{Example of network $\network_{KP-SSR}$ obtained from a state-space relaxation of $\network_{KP}$.}
\label{fig:knapsack_ssr_network}
\end{figure}

\subsection{Example on the State Space for the Capacitated Vehicle Routing Problem} \label{sec:state_space_relaxation_CVRP}

\begin{figure}[!hbt]
\center
\begin{tikzpicture}[darkstyle/.style={circle,draw,fill=gray!40,minimum size=20}, scale = 0.75, transform shape]

  \node [scale = 1.0, darkstyle, fill=white] (n1) at (0,0) {0, 0};
  \node [scale = 1.0, darkstyle, fill=white] (n2) at (-7,-2.0) {1, 1};
  \node [scale = 1.0, darkstyle, fill=white] (n3) at (-3,-2.0) {2, 2};
  \node [scale = 1.0, darkstyle, fill=white] (n4) at (3,-2.0) {2, 3};
  \node [scale = 1.0, darkstyle, fill=white] (n5) at (7,-2.0) {5, 4};
  \node [scale = 1.0, darkstyle, fill=white] (n6) at (-9,-5) {3, 2};
  \node [scale = 1.0, darkstyle, fill=white] (n7) at (-7,-5) {3, 3};
  \node [scale = 1.0, darkstyle, fill=white] (n8) at (-5,-5) {6, 4};
  \node [scale = 1.0, darkstyle, fill=white] (n9) at (-2.5,-5) {3, 1};
  \node [scale = 1.0, darkstyle, fill=white] (n10) at (1.5,-5) {4, 3};
  \node [scale = 1.0, darkstyle, fill=white] (n13) at (4.5,-5) {4, 2};
  \node [scale = 1.0, darkstyle, fill=white] (n14) at (7,-5) {6, 1};
  \node [scale = 1.0, darkstyle, fill=white] (n15) at (-6,-8.5) {5, 3};
  \node [scale = 1.0, darkstyle, fill=white] (n16) at (0,-8.5) {5, 2};
  \node [scale = 1.0, darkstyle, fill=white] (n17) at (4.5,-8.5) {5, 1};
  \node [scale = 1.4, darkstyle, fill=white] (goal) at (0,-11) {$s^*$};

  \draw[->] (n1) -- (n2) node[midway, fill=white]{-5};
  \draw[->] (n1) -- (n3) node[midway, fill=white]{-3};
  \draw[->, line width = 1.5] (n1) -- (n4) node[midway, fill=white]{1};
  \draw[->] (n1) -- (n5) node[midway, fill=white]{-2};
  \draw[->] (n2) -- (n6) node[midway, fill=white]{6};
  \draw[->] (n2) -- (n7) node[midway, fill=white]{2};
  \draw[->] (n2) -- (n8) node[midway, fill=white]{-2};
  \draw[->] (n3) -- (n9) node[midway, fill=white]{-7};
  \draw[->] (n3) -- (n10) node[midway, fill=white]{-1};
  \draw[->] (n4) -- (n9) node[midway, fill=white]{2};
  \draw[->, line width = 1.5] (n4) -- (n13) node[midway, fill=white]{-5};
  \draw[->] (n5) -- (n14) node[midway, fill=white]{-2};
  \draw[->] (n6) -- (n15) node[midway, fill=white]{-1};
  \draw[->] (n7) -- (n16) node[midway, fill=white]{-5};
  \draw[->] (n9) -- (n15) node[pos=0.6, fill=white]{2};
  \draw[->] (n10) -- (n17) node[midway, fill=white]{2};
  \draw[->] (n9) -- (n16) node[midway, fill=white]{6};
  \draw[->, line width = 1.5] (n13) -- (n17) node[midway, fill=white]{-7};

  \draw[dotted,->] (n1) -- (goal);
  \draw[dotted,->] (n2) -- (goal);
  \draw[dotted,->] (n3) -- (goal);
  \draw[dotted,->] (n4) -- (goal);
  \draw[dotted,->] (n5) -- (goal);
  \draw[dotted,->] (n6) -- (goal);
  \draw[dotted,->] (n7) -- (goal);
  \draw[dotted,->] (n8) -- (goal);
  \draw[dotted,->] (n9) -- (goal);
  \draw[dotted,->] (n10) -- (goal);
  \draw[dotted,->] (n13) -- (goal);
  \draw[dotted,->] (n14) -- (goal);
  \draw[dotted,->] (n15) -- (goal);
  \draw[dotted,->] (n16) -- (goal);
  \draw[dotted,->, line width = 1.5] (n17) -- (goal);

\end{tikzpicture}
\caption{\rev{Example of network $\network_{SSR-ESPPRC}$ obtained from a state-space relaxation of $\network_{ESPPRC}$.}}
\label{fig:ssr-espprc_network}
\end{figure}
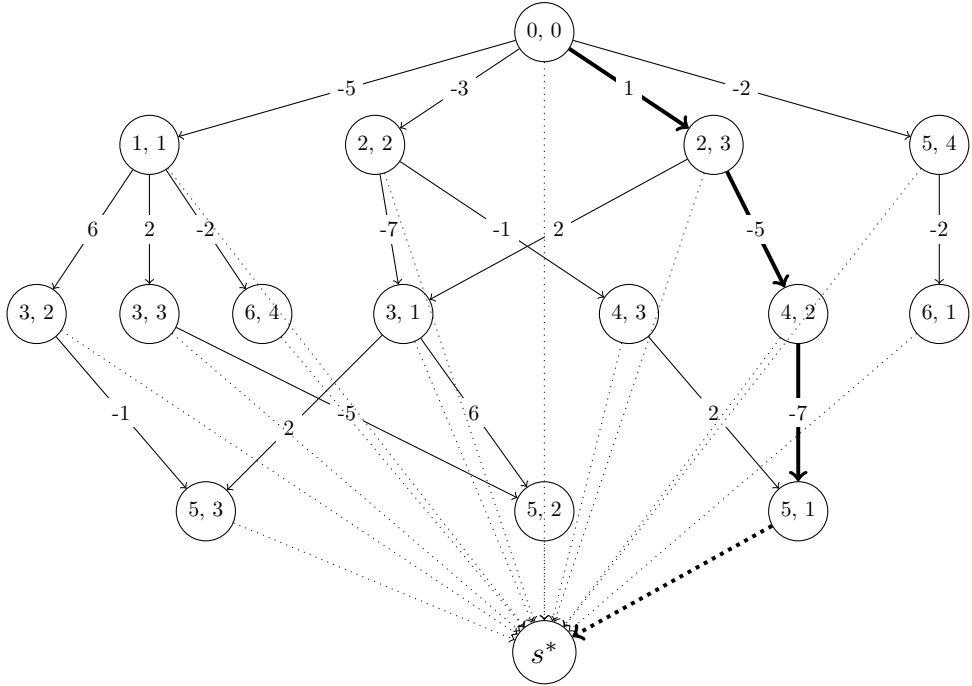

In Section \ref{sec:example_cvrp}, we presented a DW decomposition for the CVRP that resulted in the path flow model \eqref{eq:cvrp_set_partitioning_of}--\eqref{eq:cvrp_set_partitioning_domain}, which is known to have a strong linear relaxation. Then, we presented the corresponding arc flow model \eqref{eq:cvrp_exponential_arc_flow_of}--\eqref{eq:cvrp_exponential_arc_flow_domain}, which is based on the network $\network_{ESPPRC}$ for the ESPPRC and has exponential size. In the following, we present a state-space relaxation for $\network_{ESPPRC}$ that leads to a pseudo-polynomial arc flow model for the CVRP.

A pseudo-polynomial state space $\network_{ESPPRC-SSR}$ can be obtained from a state-space relaxation of $\network_{ESPPRC}$ based on the mapping function $g((S, i)) = (\sum_{j \in S} r_j, i)$, for every state $(S, i) \in \network_{ESPPRC}$. This mapping function (originally proposed by \citet{CMT81}) merges states with the same total load from the visited clients into a single state. The resulting network $\network_{ESPPRC-SSR}$, which has $O(|I|W)$ nodes and $O(|I|^2W)$ arcs, preserves the resource constraints of the ESPPRC, but it may consider non-elementary paths (clients may be visited more than once).

In the exponential arc flow model \eqref{eq:cvrp_exponential_arc_flow_of}--\eqref{eq:cvrp_exponential_arc_flow_domain} for the CVRP, the side constraints guarantee that clients are visited exactly once. Hence, by changing network $\network_{ESPPRC}$ by $\network_{ESPPRC-SSR}$ in model \eqref{eq:cvrp_exponential_arc_flow_of}--\eqref{eq:cvrp_exponential_arc_flow_domain}, the resulting model still solves the CVRP and has pseudo-polynomial size. This resulting model is often referred to as capacity-indexed formulation (see, e.g., \citet{PU14}), and it has been studied as a layered graph formulation by \citet{GLR19}.

\rev{Figure \ref{fig:ssr-espprc_network} illustrates the network $\network_{SSR-ESPPRC}$ obtained from the state-space relaxation of the network from Figure \ref{fig:ESPPRC_DP}. It can be noticed that nodes $(\{1,2\}, 1)$ and $(\{1,3\}, 1)$ have been aggregated into a single node $(3,1)$. As a result, $\network_{SSR-ESPPRC}$ has one node and one arc less than $\network_{ESPPRC}$. This is a minimal example of a reduction provided by $\network_{SSR-ESPPRC}$ that could be presented within the limits of this paper. However, due to the contrast of the exponential size of $\network_{ESPPRC}$ and the pseudo-polynomial size of $\network_{SSR-ESPPRC}$, there are many practical instances where the state-space relaxation provides a huge reduction on the size of the network.}

This example shows how state-space relaxation derives pseudo-polynomial arc flow models from a network of exponential size. This kind of modeling technique is often used in vehicle routing problems, as the combinatorial structure of such problems usually leads to exponential state spaces in DW decompositions (see, e.g., \citet{RS08}).

\section{Dual Insight} \label{sec:dual_insight}
\rev{
Research has shown that controlling the values of the dual-variables when using the primal simplex algorithm may improve computational times substantially. This happens, for instance, by modifying the model adding extra primal variables (corresponding to extra dual constraints) that must take null values at the end of the solution process. The use of this strategy may help in reducing difficulties related with instability of dual variables, primal degeneracy and long-tail effects, which are known to occur in column generation. }

\rev{
Other than the considerably smaller size, the use of more dual information is another advantage of the arc flow formulations over path flow formulations. Dual constraints (i.e., dual optimality conditions) of path flow formulations are non-negative combinations of dual constraints of arc-flow formulations. From a primal point of view, this competitive advantage can be interpreted as resulting from the possible recombination of basic variables to generate different paths in column generation algorithms (see, e.g., Sadykov and Vanderbeck \cite{SV13}).}
\subsection{On the Dual Space of Network Flow Formulations}

Many state-of-the-art algorithms to solve LP models are based on iteratively pivoting from one vertex of the LP polytope to a better neighboring vertex until an optimal solution is found. An issue that may be critical for the computational time is degeneracy. A solution is degenerate when there are basic variables with a null value; in such cases, there may be degenerate pivots that lead to the same degenerate vertex, with no change in the primal solution, nor improvement in the objective function. Hard combinatorial optimization problems are often associated with highly degenerate models, and stalling, which is a sequence of degenerate pivots, occurs in practice (see, e.g., \citet{BJS11}).
In fact, in a degenerate pivot, there is no change in the primal solution, but the new set of basic variables yields a change in the dual solution, leading to an alternative dual solution, often with a high oscillation of the values of the dual variables.

This instability often happens with master problems of DW reformulations, which may have a huge number of dual solutions associated to each primal solution. Several strategies showed that controlling the dual-variable values in the primal simplex algorithm may make the column generation procedure more efficient. For instance, \citet{MVDH99} introduced a stabilization procedure, combining a perturbation method and a penalty method, that amounts to penalizing dual variables when they lie outside a predefined box. \rev{Wentges \cite{W97} searches dual solutions in the neighborhood of the best dual solution found so far, thus reducing instability.} Other strategies include aggregating primal constraints (aggregation may change dynamically along the solution process), which enables transferring degeneracy to a complementary problem that is able to select a more central dual solution, as in \citet{EVSD05,EMDS11}. Further improvements aim at identifying a set of non-basic variables that are pivoted together into the basis, avoiding degeneracy and strictly decreasing the objective function value, by solving a problem, coined as complementary problem, in \citet{BEMS17}. For other strategies and insights on overcoming instability in combinatorial optimization algorithms, the reader is referred to, e.g., \citet{L01}.

Another strategy to deal with primal degeneracy and instability is to develop strong models with a more restricted dual space.
When comparing different LP models for the same problem with the same primal strength, the one with a tighter description of the dual space eliminates alternative dual solutions, potentially reducing degeneracy. This concept has been used in column generation algorithms by adding dual cuts that preserve all dual optimal solutions or even at least just one dual optimal solution, leading to significant speed-ups and reductions in the number of degenerate pivots (see, e.g., \citet{VC05}, \citet{LD05}, and \citet{BADVC06}).

LP solvers rely on the dual solution to prove optimality, and each dual constraint is an optimality condition. In the following, we show that arc flow formulations provide a tighter description of the dual space than their corresponding path flow formulations, with a richer description of the optimality conditions of the LP models. For instance, consider the dual of the linear relaxation of the path flow formulation \eqref{eq:general_path_formulation_of}--\eqref{eq:general_path_formulation_domain}, given by:
\begin{align}
\max ~ & \sum_{k \in {\cal C}} d_k\beta_k, & \\
\label{eq:general_path_dual_constraints0}
\text{s.t.: } & \sum_{k \in {\cal C}} \sum_{(u,v) \in \arcset(p)} \ell_{(u,v)}^k\beta_k \leq \tilde{c}_p, & \forall p \in {\cal P(N)}, \\
& \beta_k \geq 0, & \forall k \in {\cal C},
\end{align}

where dual variable $\beta_k$ correspond to the constraints \eqref{eq:general_path_formulation_constraints}. Now, consider the dual of the linear relaxation of the corresponding arc flow formulation \eqref{eq:general_arc_flow_formulation_of}--\eqref{eq:general_arc_flow_formulation_domain2}, given by:
\begin{align}
\max ~ & \sum_{k \in {\cal C}} d_k\beta_k, & \\
\label{eq:general_arc_flow_dual_constraints0}
\text{s.t.: } & \alpha_v - \alpha_u + \sum_{k \in {\cal C}} \ell_{(u,v)}^k \beta_k \leq c{(u,v)}, & \forall (u,v) \in {\arcset(\network)}, \\
\label{eq:general_arc_flow_dual_constraint1}
& \alpha_{source} - \alpha_{sink} \leq 0, & \\
& \alpha_v \geq 0, & \forall v \in \nodeset(\network), & \\
& \beta_k \geq 0, & \forall k \in {\cal C}.
\end{align}

Each dual variable $\alpha_v$ corresponds to the flow conservation of node $v \in \nodeset(\network)$, and each dual variable $\beta_k$ corresponds to the side constraint $k \in {\cal C}$. For a given path $p \in {\cal P(N)}$, by performing a non-negative linear combination of the dual constraints \eqref{eq:general_arc_flow_dual_constraints0} of every arc on $\arcset(p)$, we obtain:
\begin{align}
\label{eq:dual_linear_sum}
\sum_{(u,v) \in \arcset(p)} (\alpha_v - \alpha_u) + \sum_{(u,v) \in \arcset(p)} \sum_{k \in {\cal C}} \ell_{(u,v)}^k \beta_k \leq \sum_{(u,v) \in \arcset(p)} c{(u,v)} = \tilde{c}_p.
\end{align}

Note that, each node $v \in p$, except $v_{source}$ and $v_{sink}$, is the head of an arc $(u,v) \in \arcset(p)$ and the tail of an arc $(v,w) \in \arcset(p)$, producing in the first summation, respectively, the terms $\alpha_v$ and $-\alpha_v$ that cancel each other. Then, \eqref{eq:dual_linear_sum} can be rewritten as:
\begin{align}
(\alpha_{sink} - \alpha_{source}) + \sum_{k \in {\cal C}} \sum_{(u,v) \in \arcset(p)} \ell_{(u,v)}^k \beta_k \leq \tilde{c}_p,
\end{align}

which is equivalent to the dual constraint from \eqref{eq:general_path_dual_constraints0} for the path $p$, with an additional term $(\alpha_{sink} - \alpha_{source})$. From \eqref{eq:general_arc_flow_dual_constraint1}, we know that this term is always non-negative, implying a tighter constraint. Thus, we conclude that every dual constraint of the path flow formulation is a redundant dual constraint for the arc flow formulation, implying that the arc flow formulation provides a tighter dual space than the path flow formulation.

\rev{Generally, when the linear relaxation of either path flow or arc flow models is solved by simplex algorithms, the basis at each iteration is associated with a set of paths forming a primal-feasible solution. This set of paths is unique in the case of path flow, but not necessarily unique in the case of arc flow. If a same set of paths is considered to form the basis of an arc flow and of a path flow model, the basis of the former will be larger, due to the additional flow conservation constraints and the fact that each path is decomposed in a set of arcs in this model. But, in fact, the arc flow basis can be seen as obtained from a disaggregation of the path flow basis, which directly implies more optimality conditions. This guarantees a better description of the dual space, which, as already discussed, provides a number of practical benefits.
}

\rev{
Concluding, in practice, each simplex iteration of an arc flow model can be more expensive when compared to a path flow model (due to the larger basis), but the number of pricing iterations needed to reach proven optimality can be substantially smaller, which in many cases is a significant advantage.
}

\subsection{Example on the Cutting Stock Problem} 

We present a numerical example to compare the dual of the classical arc flow model from \citet{V99} (see Section \ref{sec:ssr_csp}) to solve the CSP and the dual of its corresponding path flow model. The dual of the linear relaxation of this arc flow model is given by:
\begin{align}
\label{eq:csp_af_dual_of}
\max & \sum_{i \in I} \beta_i, & \\
\label{eq:csp_af_dual_arcs}
\text{s.t.: } & -\alpha_v + \alpha_{v+w_i} + \beta_i \leq 0, & \forall i \in I, (v, v+w_i) \in \arcset_i(\network_{KP-SSR}), \\
\label{eq:csp_af_dual_flow}
& \alpha_{source} - \alpha_{sink} \leq 1, & \\
& \alpha_v \geq 0, & \forall v \in \nodeset(\network_{KP-SSR}), \\
\label{eq:csp_af_dual_domain}
& \beta_i \geq 0, & \forall i \in I.
\end{align}

Variables $\alpha_{v}$ are related to the flow conservation constraints of each $v \in \nodeset(\network_{KP-SSR})$, and variables $\beta_i$ are related to the demand constraint of each $i \in I$. Constraints \eqref{eq:csp_af_dual_arcs} are related to the arc variables, and constraint \eqref{eq:csp_af_dual_flow} is related to the flow variable. The dual of the linear relaxation corresponding path flow model is given by:
\begin{align}
\label{eq:csp_set_covering_dual_of}
\max & \sum_{i \in I} \beta_i, & \\
\label{eq:csp_set_covering_dual_constraints}
\text{s.t.: } & \sum_{i \in I} a_{ip} \beta_i \leq 1, & \forall p \in \pathset(\network_{KP-SSR}), \\
\label{eq:csp_set_covering_dual_domain}
& \beta_i \geq 0, & \forall i \in I.
\end{align}

The variables $\beta_i$ are related to the demand constraints of each $i \in I$, and constraints \eqref{eq:csp_set_covering_dual_constraints} are related to each path (cutting pattern) $p$ of $\network_{KP-SSR}$, where $a_{ip}$ is an integer coefficient representing the number of times item $i$ is cut from pattern $p$.

Consider again the example from Figure \ref{fig:knapsack_ssr_network}, with bin capacity $8$ and three items having $w_1 = 4$, $w_2 = 3$, $w_3 = 2$, and $d_1 = d_2 = d_3 = 1$. Tables \ref{tab:af_dual_csp} and \ref{tab:sc_dual_csp} present, respectively, model \eqref{eq:csp_af_dual_of}--\eqref{eq:csp_af_dual_domain} and \eqref{eq:csp_set_covering_dual_of}--\eqref{eq:csp_set_covering_dual_domain} for this example. In the dual of the linear relaxation of the arc flow model, $v_{source}$ and $v_{sink}$ are represented by $0$ and $S$, respectively.
In this example, the path flow model is relatively smaller than the arc flow model, which is not common for practical instances, as the former may have exponentially more primal variables. However, our goal here is only to exemplify how the dual constraints (optimality conditions) of the path flow model are redundant dual constraints for the arc flow model. This can be observed as each dual constraint of the path flow model related to a pattern can be obtained by a non-negative linear sum of the dual constraints of the arc flow model related to the arcs that form this pattern. In particular:
pattern $\{1,2\}$ can be formed by arcs $(0,4)$, $(4,7)$, $(7, S)$ and $(S,0)$;
pattern $\{1,3\}$ by arcs $(0,4)$, $(4,6)$, $(6, S)$ and $(S,0)$;
pattern $\{2,3\}$ by arcs $(0,3)$, $(3,5)$, $(5, S)$ and $(S,0)$;
pattern $\{1\}$ by arcs $(0,4)$, $(4, S)$ and $(S,0)$;
pattern $\{2\}$ by arcs $(0,3)$, $(3, S)$ and $(S,0)$;
pattern $\{3\}$ by arcs $(0,2)$, $(2, S)$ and $(S,0)$.

\begin{table}[H]
\centering
\caption{Example of a dual arc flow formulation for the CSP.}
\label{tab:af_dual_csp}
\begin{tabular}{@{}lrrrrrrrrrlrrrrl@{}}
 & \multicolumn{1}{l}{} & \multicolumn{1}{l}{$\alpha_{0}$} & \multicolumn{1}{l}{$\alpha_1$} & \multicolumn{1}{l}{$\alpha_2$} & \multicolumn{1}{l}{$\alpha_3$} & \multicolumn{1}{l}{$\alpha_4$} & \multicolumn{1}{l}{$\alpha_5$} & \multicolumn{1}{l}{$\alpha_6$} & \multicolumn{1}{l}{$\alpha_7$} & $\alpha_8$ & \multicolumn{1}{l}{$\alpha_{S}$} & \multicolumn{1}{l}{$\beta_1$} & \multicolumn{1}{l}{$\beta_2$} & \multicolumn{1}{l}{$\beta_3$} &  \\ \midrule
arc & $(0,4)$ & -1 &  &  &  & 1 &  &  &  &  &  & 1 &  &  & $\leq 0$ \\
 & $(0,3)$ & -1 &  &  & 1 &  &  &  &  &  &  &  & 1 &  & $\leq 0$ \\
 & $(4,7)$ &  &  &  &  & -1 &  &  & 1 &  &  &  & 1 &  & $\leq 0$ \\
 & $(0,2)$ & -1 &  & 1 &  &  &  &  &  &  &  &  &  & 1 & $\leq 0$ \\
 & $(3,5)$ &  &  &  & -1 &  & 1 &  &  &  &  &  &  & 1 & $\leq 0$ \\
 & $(4,6)$ &  &  &  &  & -1 &  & 1 &  &  &  &  &  & 1 & $\leq 0$ \\
 & $(0,S)$ & -1 &  &  &  &  &  &  &  &  & 1 &  &  &  & $\leq 0$ \\
 & $(2,S)$ &  &  & -1 &  &  &  &  &  &  & 1 &  &  &  & $\leq 0$ \\
 & $(3,S)$ &  &  &  & -1 &  &  &  &  &  & 1 &  &  &  & $\leq 0$ \\
 & $(4,S)$ &  &  &  &  & -1 &  &  &  &  & 1 &  &  &  & $\leq 0$ \\
 & $(5,S)$ &  &  &  &  &  & -1 &  &  &  & 1 &  &  &  & $\leq 0$ \\
 & $(6,S)$ &  &  &  &  &  &  & -1 &  &  & 1 &  &  &  & $\leq 0$ \\
 & $(7,S)$ &  &  &  &  &  &  &  & -1 &  & 1 &  &  &  & $\leq 0$ \\
flow & $(S, 0)$ & 1 &  &  &  &  &  &  &  &  & -1 &  &  &  & $\leq 1$ \\ \midrule
$\max$ & \multicolumn{1}{l}{} & \multicolumn{1}{l}{} & \multicolumn{1}{l}{} &  &  &  &  &  &  &  &  & 1 & 1 & 1 &  \\ \bottomrule
\end{tabular}
\end{table}

\begin{table}[H]
\centering
\caption{Example of a dual path formulation for the CSP.}
\label{tab:sc_dual_csp}
\begin{tabular}{@{}lrrrrl@{}}
 & \multicolumn{1}{r}{} & \multicolumn{1}{l}{$\beta_1$} & \multicolumn{1}{l}{$\beta_2$} & \multicolumn{1}{l}{$\beta_3$} &  \\ \midrule
patterns & $\{1,2\}$ & 1 & 1 &  & $\leq 1$ \\
 & $\{2,3\}$ &  & 1 & 1 & $\leq 1$ \\
 & $\{1,3\}$ & 1 &  & 1 & $\leq 1$ \\
 & $\{1\}$ & 1 &  &  & $\leq 1$ \\
 & $\{2\}$ &  & 1 &  & $\leq 1$ \\
 & $\{3\}$ &  &  & 1 & $\leq 1$ \\ \midrule
$\max$ & \multicolumn{1}{l}{} & 1 & 1 & 1 &  \\ \bottomrule
\end{tabular}
\end{table}

\section{General Solution Methods} \label{sec:methods}
As previously discussed, an advantage of arc flow formulations over their equivalent path flow formulations is that they can be much smaller and often can be solved directly by a MILP solver (which is not practical for path flow models). However, pseudo-polynomial arc flow models can still be too large, depending on the size of the parameters of an instance. In such cases, one has to rely on more sophisticated methods to solve these models. \rev{An advantage of arc flow models derived from DW decompositions is that, since their networks are based on the underlying pricing problem, it is not always necessary to load the full network in the computer memory. Instead, one can use the structure of the pricing problem to derive methods based on column generation or iterative aggregation/disaggregation to solve the problem to integer optimality while avoiding to generate the full network. Such methods can lead to an increase in practical efficiency and avoid memory overflow when solving instances associated with huge networks. }

\subsubsection*{Column Generation}
The column generation method (introduced in Section \ref{sec:dantzig_wolfe_decomposition}) \rev{solves the linear relaxation of models with a large number of variables, and} is a popular tool to solve path flow models.
The method was proposed by \citet{FF58} to solve a path flow model for a maximal multi-commodity network flow problem and three years later generalized by \citet{DW61} to solve the models resulting from DW-decompositions.
\citet{GG61, GG63} solved a path flow model for the CSP by a column generation algorithm and were the first to show the practical efficiency of the method. Since then, column generation has been the \rev{base} of several methods to solve path flow models. For references on column generation algorithms not strictly related to arc flow formulations, we refer the interested reader to the survey by \citet{LD05} and the book by \citet{DDS06}.

Column generation was applied to solve an arc flow model for the first time by \citet{V99}, who proposed a column-and-row generation algorithm. In the context of arc flow models, column-and-row generation iteratively generates arcs to enter the simplex base, while the flow conservation constraints (rows) are restricted to \rev{nodes} where there exists at least one incoming and one outgoing arc in the restricted problem. \citet{SV13} studied the column-and-row generation method and experimentally compared the solution of path flow models by column generation, with the solution of equivalent arc flow models by column-and-row generation, and observed a faster convergence of the latter, which follows the discussion in Section \ref{sec:dual_insight}. 

In column generation algorithms (and lagrangian relaxations where only flow conservation constraints are left in the master) to solve arc flow models, an LPP on the underlying network is iteratively solved, and its computational efficiency is strictly related to the network size. In this context, when the network is too large, one may rely on dynamic graph generation methods to solve the LPP, as proposed by, e.g., \citet{FH14}, to solve LPPs that arise as sub-problems from arc flow models based on large-scale time-expanded networks.

\subsubsection*{Iterative Aggregation/Disaggregation Method}

Several authors have studied state-space relaxation techniques to reduce the size of arc flow formulations (see, e.g., \citet{MAVCH11}, \citet{CHMA17}, \citet{VC12}, \citet{Boland2017}, \citet{Boland2019}, \citet{Riedler2018}).
These techniques are equivalent to applying a surrogate relaxation to the flow conservation constraints related to subsets of nodes. Theoretically, this only reduces the number of constraints. Practically speaking, after an aggregation, many arcs (variables) associated with the same decision become equivalent in the reduced network and can be merged.
From an initial relaxation, these methods use iterative techniques in which the relaxation is refined, typically by splitting nodes that have been aggregated, until the solution produced by the relaxation is feasible for the original problem, or its value is equal to a known primal bound.

\citet{MAVCH11} were the first to use these techniques in pseudo-polynomial arc flow models. The authors used two aggregations: one produces a relaxation, the other a heuristic solution.
These results were later generalized by \citet{VC12} and \citet{CHMA17}, who studied the difficulty of the different sub-problems and the performance ratio obtained by an aggregated model.
More efficient refining strategies are studied in \citet{Riedler2018}. The authors show that their path-based techniques are more effective and underline the importance of heuristic methods in the algorithm.
In early works, the elements to be aggregated were decided beforehand. In \citet{Boland2019}, the authors introduced the paradigm of \emph{dynamic discretization discovery}, in which the discretization is constructed on the fly, producing a better relaxation, by using information from the network construction process.

Such aggregation techniques have been generalized to flows in hypergraphs by \citet{benkirane:hal-02402447}
 to solve a joint rolling-stock and train selection problem for the French railway company. The authors show that even for hypergraphs, one can safely use reduced-cost filtering on aggregated variables.
Another type of relaxation is used in \citet{Nadarajah2017} to deal with problems where several constraints can be reformulated as network-flow constraints. In their model, a network is created for each constraint, and the flow conservation constraints of all networks but one are relaxed in a Lagrangian way.

\subsubsection*{Graph Reduction Methods}
An important element of efficient solution methods for arc flow models is to determine arcs that can be removed from the network without losing optimality. Removing redundant arcs is important, as a smaller network may lead to a reduction in symmetry, a tighter relaxation, and a smaller branch-and-bound tree. \rev{Many} techniques to reduce the number of arcs are problem-dependent, \rev{usually based on dominance criteria of the underlying DP}, and we discuss some of them in Section \ref{sec:applications}, under specific applications. Nonetheless, general reduction techniques have been used to enhance arc flow models, as, for instance, the reduced-cost variable-fixing method (see, e.g., \citet{PUPR10} and \citet{KLIV19}).

Given a MILP minimization model, the reduced-cost variable-fixing method performs a domain propagation based on a dual-feasible solution of the corresponding linear relaxation and an upper bound value $z_{ub}$ corresponding to an available feasible solution. The idea is that, given the objective value $z_{lb}$ of the dual-feasible solution, any integer variable with a reduced cost greater than or equal to $z_{ub} - z_{lb}$ can be removed from the model.
\citet{IDDH10} proposed an efficient reduced-cost variable fixing algorithm for eliminating arcs in network flow models that computes a bound on the reduced cost of the arcs in arc flow models based on a dual solution of the linear relaxation of its equivalent path flow model. This method was later extended by \citet{DGI20} to determine and efficiently handle pairs of sequential arcs that cannot be in the same path in an optimal solution. \rev{The impact of different dual solutions in reduced-cost variable-fixing for network flow models has been recently discussed by de Lima et al. \cite{LIM21}.}

\section{Successful Applications of Pseudo-Polynomial Arc Flow Models} \label{sec:applications}

In this section, we discuss the main applications of arc flow models, and discuss problem-dependent solution methods and reduction criteria.

\subsection{Cutting and Packing Problems}

Arc flow models have been used in a variety of cutting and packing problems, both in one and multiple dimensions.

\subsubsection*{One-dimensional Problems}
The classical arc flow model for the CSP in \citet{V99} has a node for each partial stock size, the arcs relate the items with cut positions (item arcs) or represent loss stock (loss arcs), and the combination of arcs into paths represents cutting patterns.
To solve this model to integer optimality, \citet{V99} proposed a branch-and-price algorithm based on column-and-row generation. To accelerate the column generation's convergence, the pricing problem generates paths (instead of single arcs). 
\citet{V99} proved that the arc flow model is equivalent to the path (set-covering) model by \citet{GG61, GG63}, whereas \citet{MSV18} and \citet{DI19} proved that the arc flow model is equivalent to the one-cut model by \citet{R76} and \citet{D81}. The one-cut is a pseudo-polynomial model where variables represent cutting operations on the roll.

To reduce the number of arcs, \citet{V99} constructs the graph considering that items of a single roll can always be ordered by non-increasing width. The resulting graph independently follows the state-space relaxation discussed in Section \ref{sec:ssr_csp}.
To obtain even smaller networks, \citet{CI18} proposed the meet-in-the-middle technique: each path representing a cutting pattern can be transformed into an equivalent one by left aligning the items whose left border is at the left of a given threshold parameter $t$, and right aligning the remaining items.
\rev{Recently, de Lima et al. \cite{LIM21} proposed a new way to reduce the graph in \cite{V99} by considering a maximum waste for each roll. They developed a branch-and-price framework in which the branching produces a series of small arc flow models which are solved one at a time by a general purpose MILP solver.}

\citet{CO10} proposed the DP-flow, an arc flow model for the CSP (already introduced in Section \ref{sec:example1_csp}) based on the DP network of the KP. Differently from the classical arc flow model for the CSP, which considers a node for each partial stock size, the DP-flow considers a node for each pair of items and partial stock size. The network has a level for each item, and each feasible path visits the level of each item exactly once. Although this modeling technique can substantially increase the number of nodes, it allows one to consider only proper patterns. On the other hand, the classical arc flow model is smaller, but it cannot distinguish between proper and non-proper patterns, so its relaxation can be weaker than the one of the DP-flow.

A generalization of the classical arc flow model for the CSP was proposed by \citet{BP15}. This generalization can model related problems, such as the vector BPP, the BPP with conflicts, the cardinality constrained BPP and CSP, and the graph coloring problem. The resulting network may be significantly large, but several techniques (referred to as graph compression) are proposed to reduce the network size. To break symmetry, the authors proposed a modeling technique similar to the graph construction of the DP-flow, considering an extra dimension, related to the items, on the nodes.

\citet{DI19} proposed the {\em reflect formulation}, a pseudo-polynomial model for the CSP that considers nodes and arcs only from half of the stock size. In practice, this formulation is significantly smaller than the classical arc flow model.
Besides, the reflect formulation has been proven to be as strong as the classical arc flow model without reduction criteria.
Other than the CSP, \citet{DI19} extended the reflect formulation to solve the variable-sized BPP and the BPP with item fragmentation. \citet{DDIM19} adapted the reflect formulation to solve a feasibility problem in a Benders' decomposition algorithm for the multiple knapsack problem.

\citet{AV08} presented a branch-and-price-and-cut algorithm for the multiple length CSP that solves the \citet{GG63} machine balance model, which is a path flow formulation, using the original arc flow model to generate attractive columns in a single subproblem and the variables of the arc flow model to implement a branching scheme, by expressing the branching constraints in terms of the Gilmore and Gomory model variables. The equivalence of the path and arc flow models ensures a correct transferral of dual information. Valid dual inequalities are used to stabilize and accelerate the search in the entire branch-and-bound tree.

Recently, arc flow models were proposed for the {\em skiving stock problem} (SSP), which is strongly related to the dual BPP. In the SSP, we are given a set of items, each having a length and a maximum number of copies, to be combined into the maximum number of larger items of minimum length $W$. \citet{MS16a} proposed an arc flow model for the SPP where nodes represent the length sum of combinations of items, arcs represent the positioning of the items in a combination, and paths represent a combination of items. This model, which is similar to the classical arc flow model for the CSP, cannot distinguish between proper and non-proper patterns. \citet{MS16} proposed an arc flow model that considers only proper patterns, which, similarly to the DP-flow model, considers a node for each pair of items and possible length sum, and the network has a level for each item. \citet{MDISS20} proposed two arc flow models for the SSP: one considers reversed loss arcs, which lead to a reduction of the worst-case number of nodes from $2W$ to $W$, and the other is based on the reflect model for the CSP.

\subsubsection*{Multi-dimensional Problems}

\citet{MAV10} extended the classical arc flow model to minimize the number of bins in the two-stage two-dimensional guillotine CSP.
Their model has a one-dimensional arc flow graph for the first stage cuts (which cut the bins to obtain strips), and a one-dimensional arc flow graph for each possible strip size to determine the second stage cuts (which produce the items from the strips, possibly admitting a final trim loss cut).
To reduce symmetry, the flow in the second stage graph of a given strip size is equal to the sum of the flows for that strip size in the first stage graph. Note that these strips may belong to the same bin or different bins.
Procedures to reduce the size of the graph and a new family of cutting planes based on the height of the items were proposed.
The arc flow model in \cite{MAV10} was later adapted by \citet{M15} to solve the two-stage two-dimensional guillotine strip packing problem.

\citet{NDIS17} proposed an arc flow model, similar to the one by \citet{MAV10}, to solve a three-stage two-dimensional strip packing problem where a limit is imposed on the number of shelves and setup times between items must be taken into account.
\citet{DIM17} adapted the classical arc flow model for the CSP to solve the one-dimensional contiguous bin packing problem, which often appears as a sub-problem in two-dimensional cutting and packing solution methods. The authors used it to solve the sub-problem of a Benders' decomposition algorithm for the two-dimensional strip packing with item rotations and for the pallet loading problem.

\citet{CSVV18} solved the four-stage two-dimensional guillotine bounded knapsack problem with a network flow model based on a directed acyclic hypergraph. They compared the efficiency of several algorithms based on this representation, including a MILP model and an iterative state-space relaxation based on the corresponding DP.

\subsection{Scheduling Problems}

In the scheduling field, both time-indexed (see, e.g., \citet{SW92}) and path flow (see, e.g., \citet{VHS00}) formulations have been widely used to solve a variety of optimization problems.
In a closely related context, polynomial arc flow models have been proposed more than three decades ago by \citet{EM87} to solve a lot-sizing problem. However, pseudo-polynomial arc flow formulations have only been adopted for scheduling problems during the last decade.
It is worth mentioning that time-indexed formulations can be used to derive arc flow models of the same strength, which are associated to a sparser constraint matrix. The equivalence between time-indexed and arc flow formulations has been shown in different contexts starting from \citet{V02}, who proved such equivalence by a unimodular transformation, in the context of the CSP.

\citet{PUPR10} developed a branch-and-cut-and-price algorithm for the problem of minimizing weighted tardiness on identical parallel machines (denoted as $P||\sum w_j T_j$ in the three-field classification of \citet{GLLK79}). Their algorithm is based on an arc-time-indexed formulation and is improved with a number of combinatorial techniques, including variable fixing by reduced costs, extended capacity cuts, dual stabilization, and the direct solution of the formulation by a MILP solver if the fixing procedure had consistently reduced the number of variables. The method in \cite{PUPR10} was later extended by \citet{BSSU20}, who proposed a branch-and-cut-and-price algorithm to solve a path flow formulation for parallel machine scheduling in which the branching is based on the variables of the arc flow model.

\citet{LRS11} developed a branch-and-price algorithm for the job shop problem with a general min-sum objective function. Their algorithm is based on the solution of an arc flow model in which a path has to be chosen for each job (which is composed of multiple operations). The model was strengthened by clique inequalities.

\citet{RBMHW13} presented a comprehensive list of mathematical models for scheduling jobs on a single machine by minimizing weighted earliness and tardiness. The problem, denoted as $1||\sum \alpha_j E_j + \sum \beta_j T_j$, is relevant in the context of just-in-time production. In computational tests on random instances, the arc flow model achieved the lowest optimality gaps.

\citet{MS18} presented a direct extension of the arc flow formulation by \citet{V99} to the problem of scheduling jobs on identical parallel machines with the objective of minimizing the makespan ($P||C_{\max}$).

\citet{KDI19} considered again the problem of scheduling jobs on identical parallel machines, but focused on the minimization of the weighted sum of the completion times ($P||\sum w_j C_j$). They presented an arc flow model, and then enhanced it by grouping jobs having the same weight and processing time and creating time windows for each group by considering job priorities. A computational comparison showed that the enhanced arc flow performed very well compared to other time-indexed, convex integer quadratic programming, and path flow models.

\citet{KIL19} extended the work in \cite{KDI19} to deal with the case of family setup times ($P|s_i| \sum w_j C_j$). The authors proposed three different arc flow models. In the most efficient one, the network is divided into a set of layers, one layer per family. Arcs within the same layer considered only the processing time of a job, whereas arcs connecting two layers considered both setup and processing times. Computational results showed that setup times worsen the linear relaxation value of the arc flow models, which outperformed a path flow model only on a handful of instances. 

A further generalization of \cite{KDI19} was provided in \citet{KDFI20}, who considered the case of release dates ($P|r_j|\sum w_j C_j$) and obtained stricter time windows for the jobs. The resulting arc flow model obtained better results than a branch-and-price based on the path flow formulation on instances having jobs of small and moderate processing time.

We also mention that interesting integrations between scheduling and other combinatorial problems have been tackled in the literature. 
\citet{CS15} considered a joint assignment, scheduling, and routing problem arising in home care optimization. An arc flow model was used for the generation of patterns, which represent feasible combinations of the three decision levels of the problem.
\citet{CDYC19} proposed an arc flow formulation for a rotation assignment and scheduling problem arising in the context of clinical rotations. They demonstrate that the network model was computationally superior to a classical MILP model on a real-world set of instances.

\citet{BAMV16} compared two formulations for a combined cutting stock and scheduling problem: a compact formulation strengthened with knapsack inequalities and an arc flow formulation. Using a revised version of the latter formulation based on aggregated time periods, they derived a heuristic solution procedure that proved to be effective for the solution of medium size instances.
\citet{RABV16} proposed and analyzed an arc flow formulation for a combined cutting stock and scheduling problem on parallel machines. Different strategies to simplify the formulation by reducing the number of arcs were presented. 

\citet{ROM20} solved scheduling problems with batch processing machines. In such problems, the jobs are grouped into batches to be scheduled in machines, where the batches have a capacity to be respected by the size of its grouped items. The author proposed arc flow models where nodes represent a discretization of the capacity of a batch and arcs represent either a scheduled job in a batch or an unused capacity.

\subsection{Routing Problems}

In routing problems, the input is usually based on a network, where clients and depots are given as nodes and arcs are related to transportation between the nodes (see, e.g., \citet{TV14}). For such problems, arc flow models based on the input network lead to compact models that are small but may have weak relaxations and too much symmetry to be solved in practice. Examples of such compact models are the classical MTZ model by \citet{MTZ60} for the traveling salesman problem (TSP) and the three-index formulation (see Section \ref{sec:example_cvrp}) for the CVRP. In many routing problems, the ability to avoid non-elementary routes is an important aspect to determine the strength of a relaxation. For this reason, the best solution methods for many variants, especially the ones with multiple vehicles, are usually based on path flow formulations (see, e.g., \citet{PU14} and \citet{PSUV19}), as they can handle non-elementary routes more easily than arc flow formulations. However, pseudo-polynomial arc flow formulations have still been used to solve open problems of a number of variants, and to enhance state-of-the-art methods based on path flow formulations.

Pseudo-polynomial arc flow formulations for vehicle routing problems were first introduced by \citet{GGMPP07} and \citet{PPU08}, which, inspired by the early work by \citet{PQ78} for single machine scheduling with minimum tardiness, proposed the capacity-indexed formulation for the CVRP. The capacity-indexed formulation, which follows the network from the state-space relaxation in Section \ref{sec:state_space_relaxation_CVRP}, has a node $(i,q)$ for each pair of client (or depot) $i$ and partial capacity $q$, where the partial capacities group the nodes into levels (layers). Then, each path arriving at a node $(i,q)$ represents a route that finishes in client $i$ and has total load $q$.

The capacity-indexed formulation is based on a network where paths may be associated with non-elementary routes (which weakens its linear relaxation), and may be too large to be solved in practice (see, e.g., \citet{PPU08}). The main interest in this formulation is that its variables can be used to define cutting planes (e.g., the extended capacity cuts) to strengthen path flow formulations for vehicle routing problems (see, e.g., \citet{PU14}). 
According to \citet{U11}, cutting planes based on the capacity-indexed formulation have been successfully used in other vehicle routing problems, and even in parallel machine scheduling problems.

\citet{MAVCH11} proposed an iterative aggregation/disaggregation algorithm to solve an arc flow formulation for the vehicle routing problem with time windows and multiple routes.
In the underlying network, each node corresponds to a time instant, and each arc corresponds to a possible subtour.
Aggregation based on rounding procedures was used to make the routes with non-integer travel times fit the model's graph. Whenever an infeasible solution was found, the nodes involved in the infeasibility were disaggregated, and the resulting model was solved again.

\citet{BAM17} proposed an arc flow formulation for the multi-trip inventory routing problem where vehicles can perform more than a single route per time period. Following \citet{MAVCH11}, nodes and arcs of the underlying network correspond to, respectively, time instants and routes that may be assigned to a vehicle.

Algorithms based on iterative aggregation/disaggregation were later proposed to solve time-expanded arc flow formulations for the TSP with time windows (TSPTW). The network of time-expanded formulations has a node $(i,t)$ for each client $i$ and time instant $t$. Each path from the source to a node $(i,t)$ represents a path arriving at client $i$ at time $t$.
\citet{BHSD17} and \citet{Riedler2018} proposed iterative aggregation/disaggregation algorithms to solve a time-expanded model for the TSPTW, starting with a reduced network that is sufficient to produce a bound for the problem and is iteratively refined to find an optimal solution. The method in \citet{BHSD17} was later extended by \citet{VHBS19} to solve the generalization where travel times are time-dependent. This kind of approach was also used by \citet{Boland2017} to solve a continuous-time service network design problem without the need to approximate the solution by a discretization of the time-horizon.

\subsection{Miscellaneous}

Earth Observation Satellite scheduling requires to determine the pictures to be taken by a set of satellites in a given time period, so as to satisfy side constraints and optimize an objective function. Pseudo-polynomial arc flow models have been proposed to solve such problems by \citet{GM03} and \citet{WWXP20}. In these models, the nodes are related to a discretization of the time horizon, and the arcs are related to the decisions on the pictures to be taken.

\citet{KLIV19} solved the dynamic berth allocation problem, which aims at allocating vessels into quays that are divided into berths, while optimizing an objective function based on the service time of each vessel. The authors proposed an arc flow model where nodes are associated with time instants and arcs are related to vessels serving. Problem-dependent reduction criteria and a reduced-cost variable-fixing algorithm were proposed to improve the solution time.

In the capacitated $p$-center problem, a set of customers must be attributed to capacitated facility locations, minimizing the maximum distance between each client and its facility. To solve this problem, \citet{KIV19} proposed an arc flow model where the underlying graph has a component for each location, in which nodes correspond to partial filling of the facility capacity, and arcs correspond to customers and unused capacity.

\citet{RAV20} described an arc flow formulation for the multi-trip production, inventory, distribution, and routing problem with time windows. Nodes and arcs represent time instants and vehicle routes, respectively. 

Train timetabling problems require to determine a periodic timetable for a set of trains that satisfies operational constraints and optimizes an objective function (see, e.g., \citet{CT12}). The input for such problems is based on a graph where nodes represent stations and arcs represent tracks. Train timetabling problems have been successfully solved by arc flow models based on time-expanded networks by \citet{CFT02} and \citet{FH14}.

Recently, \citet{V20} proposed an arc flow model to solve the graph coloring problem, which asks to partition a graph into the minimum number of independent sets. The proposed arc flow model is based on decision diagrams, which are closely related to DP (see, e.g., \citet{H13}) and result in acyclic graphs. The resulting network in \cite{V20} is equivalent to a DP network to solve the maximum independent set problem, which is usually the pricing problem of path (set-partitioning) formulations for the graph coloring problem (see, e.g., \citet{GM12}). An iterative refinement method is proposed by the authors to deal with the exponential size of the model.

\section{Conclusion and Future Research Directions} \label{sec:conclusion}

In this survey, we reviewed over one hundred papers related to arc flow formulations. Many of these papers present arc flow models having pseudo-polynomial size and a strong linear relaxation. The number of applications of pseudo-polynomial arc flow formulations has grown considerably since the work by \citet{V99}, making them a valid alternative to path flow formulations.
For many combinatorial optimization problems, path flow formulations are still the best alternative, because they can embed difficult constraints in the subproblem solved to build the paths. However, arc flow formulations have many positive aspects, among which we highlight that: 
they are powerful modeling tools that allow one to model complex issues from real systems;
pseudo-polynomial arc flow models are often related to DW decompositions, providing strong linear relaxations;
differently from path flow formulations, they provide models that usually have a number of variables which allows practical solutions directly by general MILP solvers, avoiding complex implementations;
compared to equivalent path flow formulations, they have a richer description of the dual space, leading to a faster convergence of simplex-based methods;
pseudo-polynomial arc flow models can be derived from state-space relaxations from the underlying DP network from path flow models.

This survey is a contribution to a systematic study of arc flow formulations, but we would like to point out that there are many open questions and research lines to be pursued. Some of them are the following.

One of the advantages of arc flow formulations is that they provide a tighter description of the dual feasible space keeping the primal strength. Besides the methodologies presented in this paper, are there any general hints on how to do both primal and dual strengthening in arc flow models? For example, in extended formulations, using new sets of variables may strengthen the primal model and, as new (primal) variables are dual cuts, there is also a richer description of the dual feasible space.

Models presented in this survey explore solution spaces that are convex combinations of flows, each corresponding to an $s-t$ path, which is a sequence of arcs and vertices. Are there other types of structures (e.g., involving sequences of operations) that also lead to models with strong bounds? In fact, there are pseudo-polynomial models that, instead of using a sequence of arcs to form a path (which is an extremal solution), use a sequence of operations to form a solution that is extremal. Examples are \citet{D81} for the CSP, and \citet{SAV10} and \citet{FMT16} for the two-dimensional CSP with guillotine cuts. In these cases, one-cut operations are/have to be combined to form a cutting pattern (the extremal solution).

There are pseudo-polynomial models that do not provide LP lower bounds as strong as those of column generation (e.g., position indexed models for two-dimensional non-guillotine CSP). Is there any structural (extremal) property, in these cases, that can be explored and may lead to models with stronger bounds?

Several solution methods for large-scale arc flow models, like the ones presented in Section \ref{sec:methods}, are general and can be applied to any arc flow model based on DP. A \rev{software/library containing such general tools to solve general large-scale arc flow models} would be an interesting contribution to the optimization and operations research community.

\section*{Acknowledgement}

We thank three anonymous referees for their careful reviews. Their comments helped improve the clarity of the presentation and the overall quality of the paper.
The first and fourth authors have been supported by FAPESP - Funda{\c{c}}{\~a}o de Amparo {\`a} Pesquisa do Estado de S{\~a}o Paulo (under grant numbers 2017/11831-1 and 2019/12728-5). The second and the fifth authors have been supported by FCT - Funda{\c{c}}{\~a}o para a Ci{\^e}ncia e Tecnologia within the R\&D Units Project Scope: UIDB/00319/2020. 

\bibliographystyle{myplainnat}

\end{document}